\documentclass[11pt]{amsart}
\usepackage{amsmath}
\usepackage{amsxtra}
\usepackage{amscd}
\usepackage{amsthm}
\usepackage{amsfonts}
\usepackage{amssymb}
\usepackage{eucal}

\newtheorem{lem}{Lemma}
\newtheorem{prop}{Proposition}

\newtheorem{conj}{Conjecture}

\theoremstyle{definition}

\theoremstyle{definition}
\newtheorem{thm}{Theorem}

\theoremstyle{remark}
\newtheorem{rem}{Remark}

\numberwithin{equation}{section}
\numberwithin{conj}{section}
\numberwithin{thm}{section}
\numberwithin{prop}{section}
\numberwithin{lem}{section}
\numberwithin{rem}{section}
\newcommand{\thmref}[1]{Theorem~\ref{#1}}
\newcommand{\secref}[1]{Sect.~\ref{#1}}
\newcommand{\lemref}[1]{Lemma~\ref{#1}}
\newcommand{\propref}[1]{Proposition~\ref{#1}}

\newcommand{\conjref}[1]{Conjecture~\ref{#1}}
\newcommand{\remref}[1]{Remark~\ref{#1}}
\newcommand{\nc}{\newcommand}
\nc{\on}{\operatorname}
\nc{\Z}{{\mathbb Z}}
\nc{\C}{{\mathbb C}}
\nc{\F}{{\mathbb F}_q}
\nc{\Ql}{\ol{{\mathbb Q}}_\ell}
\nc{\cond}{|\,}
\nc{\pone}{{\mathbb P}^1}
\nc{\pa}{\partial}
\nc{\arr}{\rightarrow}
\nc{\al}{\alpha}
\nc{\ri}{\rangle}
\nc{\lef}{\langle}
\nc{\W}{\Phi}
\nc{\p}{{\mathcal p}}
\nc{\la}{\lambda}
\nc{\ep}{\epsilon}
\nc{\V}{{\mathcal V}}
\nc{\De}{\Delta}
\nc{\G}{{{\mathcal G}r}}
\nc{\Li}{{\mathcal L}}
\nc{\Pro}{{\mathbb P}}
\nc{\z}{{\mathcal Z}}
\nc{\sw}{{\mathfrak s}{\mathfrak l}}
\nc{\wh}{\widehat}
\nc{\zz}{{\mathcal Z}}
\nc{\wt}{\widetilde}
\nc{\hi}{{\mathcal H}}
\nc{\M}{{\mathcal M}}
\nc{\larr}{\longrightarrow}
\nc{\La}{\Lambda}
\nc{\us}{\underset}
\nc{\PL}{{}^L P}
\nc{\GL}{{}^L G}
\nc{\TL}{{}^L T}
\nc{\A}{{\mathbb A}}
\nc{\K}{{\mathcal K}}
\nc{\OO}{{\mathcal O}}
\nc{\Q}{{\mathcal Q}}
\nc{\ol}{\overline}
\nc{\Coh}{{{\mathcal C}oh}}
\nc{\Cohn}{\Coh_n}
\nc{\f}{{\mathcal F}}
\nc{\si}{_E}
\nc{\ga}{\gamma}
\nc{\Fq}{{\mathbb F}_q}
\nc{\s}{{\mathbf s}}
\nc{\Ga}{{\mathbb G}_a}
\nc{\gen}{\on{gen}}
\nc{\GK}{G(\K)}
\nc{\GO}{\GG(\OO)}
\nc{\NK}{\NN(\K)}
\nc{\NO}{\NN(\OO)}
\nc{\cf}{{\mathcal F}}
\nc{\HH}{{\mathcal Four}^{-1}}
\nc{\Gaf}{{\mathbb G}_{a,\Fq}}
\nc{\Gafb}{{\mathbb G}_{a,\ol{\mathbb F}_q}}
\nc{\vf}{\varphi}
\nc{\GG}{{\mathcal G}}
\nc{\NN}{{\mathcal N}}
\nc{\loc}{_{\on{loc}}}
\nc{\Th}{\Theta}
\nc{\g}{^L\!{\mathfrak g}}
\nc{\bi}{\bibitem}

\nc{\Fqb}{\overline{{\mathbb F}}_q}
\nc{\Fqp}{{\mathbb F}_{q'}}

\begin{document}

\title[Geometric realization of Whittaker functions]{Geometric realization
of Whittaker functions and the Langlands conjecture}

\author{E. Frenkel}

\address{Department of Mathematics, Harvard University, Cambridge, MA
02138, USA}

\author{D. Gaitsgory}

\address{School of Mathematics, Institute for Advanced Study, Princeton, NJ
08540, USA}

\author{D. Kazhdan}

\address{Department of Mathematics, Harvard University, Cambridge, MA
02138, USA}

\author{K. Vilonen}

\address{Department of Mathematics, Brandeis University, Waltham, MA 02254,
USA}

\date{March 1997; Revised: November 1997}

\begin{abstract} We prove the equivalence of two conjectural constructions
of unramified cuspidal automorphic functions on the adelic group $GL_n(\A)$
associated to an irreducible $\ell$--adic local system of rank $n$ on an
algebraic curve $X$ over a finite field. The existence of such a function
is predicted by the Langlands conjecture.

The first construction, which was proposed by Shalika \cite{Sha} and
Piatetski-Shapiro \cite{PS} following Weil \cite{W} and Jacquet-Lang\-lands
\cite{JL} ($n=2$), is based on considering the Whittaker function. The
second construction, which was proposed recently by Laumon \cite{La2}
following Drinfeld \cite{Dr} ($n=2$) and Deligne ($n=1$), is geometric: the
automorphic function is obtained via Grothendieck's ``faisceaux-fonctions''
correspondence from a complex of sheaves on an algebraic stack.

Our proof of their equivalence is based on a local result about the
spherical Hecke algebra, which we prove for an arbitrary reductive
group. We also discuss a geometric interpretation of this result.
\end{abstract}

\maketitle

\section{Introduction}

\subsection{} Let $X$ be a smooth, complete, geometrically connected curve
over $\Fq$. Denote by $F$ the field of rational functions on $X$, by
${\mathbb A}$ the ring of adeles of $F$, and by $\on{Gal}(\ol{F}/F)$
the Galois group of $F$.

The present paper may be considered as a step towards understanding
the geometric Langlands correspondence between $n$--dimensional
$\ell$--adic representations of $\on{Gal}(\ol{F}/F)$ and automorphic
forms on the group $GL_n({\mathbb A})$. We follow the approach
initiated by V.~Drinfeld \cite{Dr}, who applied the theory of
$\ell$--adic sheaves to establish this correspondence in the case of
$GL_2$.

The starting point of Drinfeld's approach is the observation that an
unramified automorphic form on the group $GL_n({\mathbb A})$ can be viewed
as a function on the set $M_n$ of isomorphism classes of rank $n$ bundles
on the curve $X$. The set $M_n$ is the set of $\Fq$--points of ${\mathcal
M}_n$, the algebraic stack of rank $n$ bundles on $X$. One may hope to
construct the automorphic form associated to a Galois representation as a
function corresponding to an $\ell$--adic perverse sheaf on ${\mathcal
M}_n$. This is essentially what Drinfeld did in \cite{Dr} in the case of
$GL_2$. In abelian class field theory (the case of $GL_1$) this  was done
previously by P.~Deligne (see \cite{La1}).

\subsection{} Let $M'_n$ denote the set of isomorphism classes of pairs
$\{L,s\}$, where $L\in M_n$ is a rank $n$ bundle on $X$ and $s$ is a
regular non-zero section of $L$. Using a well-known construction due
to Weil \cite{W} and Jacquet-Langlands \cite{JL} for $n=2$, and
Shalika \cite{Sha} and Piatetski-Shapiro \cite{PS} for general $n$,
one can associate to an unramified $n$--dimensional representation
$\sigma$ of $\on{Gal}(\ol{F}/F)$, a function $f'_\sigma$ on
$M'_n$. The construction of $f'_\sigma$ is obtained from the
Whittaker function $W_{\sigma}$, a function canonically attached to
$\sigma$. The Langlands conjecture predicts that when $\sigma$ is
geometrically irreducible, the function $f'_\sigma$ is constant along
the fibers of the projection $p: M'_n\to M_n$. In other words,
conjecturally, $f'_\sigma$ is the pull-back of a function $f_{\sigma}$
on $M_n$; the function $f_{\sigma}$ is then the automorphic function
corresponding to $\sigma$.

\subsection{} Let now $\M'_n$ be the moduli stack of pairs $\{L,s\}$, where
$L$ is a rank
$n$ bundle on $X$ and $s$ is a regular non-zero section of $L$. We have:
$M_n=\M_n(\Fq),M'_n=\M'_n(\Fq)$. Each Galois representation $\sigma$ gives
rise to an $\ell$--adic local system $E$ on $X$ of rank $n$. Drinfeld's
idea, developed further by G.~Laumon \cite{La1}, can be interpreted as
follows.

Suppose there exists an irreducible perverse sheaf ${\mathcal S}'\si$ on
$\M'_n$, with the property that the function $S'\si$ associated to
${\mathcal S}'\si$ on $M'_n=\M'_n(\Fq)$ equals $f'_\sigma$. Then showing
that the function $f'_\sigma$ is constant along the fibers of the
projection $p: M'_n\to M_n$ becomes a geometric problem of proving that
${\mathcal S}'\si$ descends to a perverse sheaf ${\mathcal S}\si$ on $\M_n$.

In \cite{Dr}, Drinfeld constructed such a sheaf ${\mathcal S}'\si$ in
the case of $GL_2$. He started with a geometric realization of the
Whittaker function as a perverse sheaf on the symmetric power of the
curve $X$. Then he defined the sheaf ${\mathcal S}'\si$ using a
geometric version of the Weil-Jacquet-Langlands construction. Drinfeld
showed that the sheaf ${\mathcal S}'\si$ is locally constant along the
fibers of $p$. Since the fibers of $p$ are projective spaces, hence
simply-connected, this implies that ${\mathcal S}'\si$ is constant
along the fibers of $p$. One may hope to use a similar argument in the
case of $GL_n$.

\subsection{} In order to construct ${\mathcal S}'\si$ in general, Laumon
\cite{La1} defined a sheaf ${\mathcal L}\si$, which he considered as a
geometric analogue of the Whittaker function $W_\sigma$. However, the
function $L\si$ corresponding to ${\mathcal L}\si$ and the Whittaker
function $W_\sigma$ are defined on different sets and their values are
different, see \cite{La1} and Sect. 3.5 below. Using the sheaf
${\mathcal L}\si$, Laumon \cite{La2} constructed a candidate for the
sheaf ${\mathcal S}'\si$ on $\M'_n$ (this construction was
independently found by one of us; D.G., unpublished). In order to
justify this construction, one has to prove that the function $S'\si$
on $M'_n$, corresponding to the sheaf ${\mathcal S}'\si$, coincides
with the function $f'_\sigma$. This equality was conjectured by Laumon
in \cite{La2} (Conjecture 3.2), and its proof is one of the main goals
of this paper.

To prove the equality $S'\si=f'_\sigma$, we reduce it to a local
statement (see \thmref{local}) which we make for an arbitrary
reductive group. We prove \thmref{local} using the Casselman-Shalika
formula for the Whittaker function \cite{CS} (actually, \thmref{local}
is equivalent to the Casselman-Shalika formula). \thmref{local} can be
translated into a geometric statement about  intersection
cohomology sheaves on the affine Grassmannian (see \conjref{second}).

\subsection{} One essential difference between Laumon's approach and our
approach is in interpretation of the sheaf ${\mathcal L}\si$ and the
function $L\si$ associated to this sheaf. Laumon interprets the local
factors of $L\si$ via the Springer sheaves and the Kostka-Foulkes
polynomials (see \remref{3}). We interpret the local factors of $L\si$ via
the perverse sheaves on the affine Grassmannian and the Hecke algebra (see
Sect.~4.2). Our interpretation, which was inspired by \cite{Lu1}, allows us
to gain more insight into Laumon's construction. In particular, it helps to
explain the apparent discrepancy between $L\si$ and $W_\sigma$: it turns
out that $L\si$ is  related to $W_\sigma$ by a Fourier transform. Using
this result, we demonstrate that the outputs of the two constructions --
$S'\si$ and $f'_\sigma$ -- coincide.

\subsection{} Let us now briefly describe the contents of the paper:

In Sect.~2 we review some background material concerning the Langlands
conjecture and the classical construction of the function $f'_\sigma$
together with its geometric interpretation. We follow closely Sect.~1
of \cite{La1}.

In Sect.~3 we describe the construction of the sheaf ${\mathcal S}'\si$ on
$\M'_n$ and the function $S'\si$ on $M'_n$. We state the main conjecture
(\conjref{princ}) about the geometric Langlands correspondence for $GL_n$
and our main result (\thmref{prin}).

In Sect.~4 we give an adelic interpretation of the construction of the
function $S'\si$ and reduce \thmref{prin} to a local statement,
\propref{Fplus}.

In Sect.~5 we prove \thmref{local} for an arbitrary reductive group and
derive from it \propref{Fplus}.

In Sect.~6 we interpret the results of Sect.~5 from the point of view of
the spherical functions.

In Sect.~7 we give a geometric interpretation of \thmref{local} and discuss
a possible generalization of Laumon's construction to other groups.

\subsection{} In this paper we work with algebraic stacks in the smooth
topology in the sense of \cite{La:s}. All stacks that we consider have
locally the form ${\mathcal Y}/G$ where ${\mathcal Y}$ is a scheme and $G$
is an algebraic group acting on it. We will use the following notion of
perverse sheaves on such algebraic stacks: for ${\mathcal V}={\mathcal
Y}/G$, a perverse sheaf on ${\mathcal V}$ is just a $G$--equivariant
perverse sheaf on ${\mathcal Y}$, approprietly shifted.

Throughout this paper, for an $\Fq$--scheme (resp., for an
$\Fq$--stack) ${\mathcal V}$ and for an algebra $R$ over $\Fq$,
${\mathcal V}(R)$ will denote the set of $R$--points of ${\mathcal V}$
(resp., the set of isomorphism classes of $R$--points of ${\mathcal
V}$). In most cases, schemes and stacks are denoted by script letters
and their sets of ${\mathbb F}_q$--points are denoted by the
corresponding roman letters (e.g., ${\mathcal V}$ and $V$). We use the
same notation for a morphism of stacks and for the corresponding map
of sets.

If ${\mathcal S}$ is a sheaf or a complex of sheaves on a stack ${\mathcal
V}$, then the corresponding set $V$ is endowed with a function of
``alternating sum of traces of the Frobenius'' (as in \cite{De})  which we
denote by the corresponding roman letter $S$ (we assume that a square root
of $q$ in $\Ql$ is fixed throughout the paper).

If ${\mathcal V}$ is a stack over ${\mathbb F}_q$, the set $V$ is endowed
with a canonical measure $\mu$, which in the case when ${\mathcal V}$ is a
scheme has the property $\mu(v)=1$, $\forall v\in V$. For example, if
${{\mathcal G}}$ is a group, $\mu((pt/{{\mathcal G}})({\mathbb
F}_q))=|G|^{-1}$.

\subsection{Acknowledgments} We express our gratitude to J.~Bernstein and
I.~Mirkovi\'c for valuable discussions and to B.~Gross for useful
remarks concerning the Whittaker models. We thank A.~Beilinson and
V.~Drinfeld for sharing with us their ideas and unpublished results
about the affine Grassmannian, which we used in Sect.~7.2.

We are indebted to the referee for  valuable comments and suggestions.

The research of E.~Frenkel was supported by grants from the Packard
and Sloan Foundations, and by the NSF grants DMS 9501414 and DMS
9304580. The research of D.~Gaitsgory was supported by the NSF grant
DMS 9304580. D.~Kazhdan was supported by the NSF grant DMS
9622742. K.~Vilonen was supported by the NSF grant DMS 9504299 and by
the NSA grant MDA 90495H103.

\section{Background and the Shalika-Piatetski-Shapiro construction}

\subsection{Langlands conjecture} Let $k={\mathbb F}_q$ be a finite field,
and  let $X$ be a smooth complete geometrically connected curve over
$k$. Denote by $F$ the field of rational functions on $X$. For each
closed point $x$ of $X$, denote by $\K_x$ the completion of $F$ at
$x$, by $\OO_x$ the ring of integers of $\K_x$, and by $\pi_x$ a
generator of the maximal ideal of $\OO_x$. Let $k_x$ be the residue
field $\OO_x/\pi_x \OO_x$, and $q_x = q^{\deg x}$ be its
cardinality. We denote by $\A = \prod'_{x \in |X|} \K_x$ the ring of
adeles of $F$ and by $\OO = \prod'_{x \in |X|} \OO_x$ its maximal
compact subring.

Consider the set ${\mathfrak G}_n$ of unramified and geometrically
irreducible $\ell$--adic representations of the Galois group
$\on{Gal}(\overline{F}/F)$ in $GL_n(\Ql)$, where $\ell$ is relatively
prime to $q$, as in \cite{La1}, (1.1).

Let ${\mathfrak A}_n$ be the set of cuspidal unramified automorphic
functions on the group $GL_n(\A)$ -- these are cuspidal functions on
the set $GL_n(F)\backslash GL_n(\A)/GL_n(\OO)$, which are
eigenfunctions of the Hecke operators. Recall that for each $x \in
|X|$ and $i=1,\ldots,n$, one defines the Hecke operator $T^i_x$ by the
formula: $$\left( T^i_x \cdot f \right)(g) = \int_{M^i_n(\OO_x)} f(gh)
dh,$$ where $$M^i_n(\OO_x) = GL_n(\OO_x) \cdot D_x^i \cdot GL_n(\OO_x)
\subset GL_n(\K_x) \subset GL_n(\A),$$ $D_x^i$ is the diagonal matrix
whose first $i$ entries equal $\pi_x$, and the remaining $n-i$ entries
equal $1$, and $dh$ stands for the Haar measure on $GL_n(\K_x)$
normalized so that $GL_n(\OO_x)$ has measure $1$.

Let $B$ be the Borel subgroup of upper triangular matrices and $T N$
be its standard Levi decomposition. Cuspidality condition means that
for each proper parabolic subgroup of $GL_n$, whose unipotent radical
$V$ is contained in the upper unipotent subgroup $N$, 
$$\int_{V(F) \backslash V(\A)} f(vg) dv = 0, \quad \quad \forall g \in
GL_n(\A).$$

\begin{conj}    \label{Lan} For each $\sigma \in {\mathfrak G}_n$, there
exists a unique (up to a non-zero constant multiple) function
$f_\sigma \in {\mathfrak A}_n$, such that for any $x \in |X|$
$$T^i_x \cdot f_\sigma = q_x^{-i(i-1)/2} \on{Tr}(\Lambda^i
\sigma(\on{Fr}_x)) f_\sigma, \quad \quad i=1,\ldots,n,$$ where $T^i_x$ is
the $i$th Hecke operator, and $\on{Fr}_x \in \on{Gal}(\overline{F}/F)$ is
 the geometric Frobenius element.
\end{conj}

Let $P_1 \subset GL_n$ be the subgroup
\begin{equation}    \label{pp}
\left\{ \begin{pmatrix} g & h \\ 0 & 1
\end{pmatrix} \cond g \in GL_{n-1} \right\}.
\end{equation} Following Shalika and Piatetski-Shapiro, we construct a
function
$f'_\sigma$ on the double-quotient
$$P_1(F)\backslash GL_n(\A)/GL_n(\OO),$$ which is cuspidal and which
satisfies the Hecke eigenfunction property:
\begin{equation}    \label{hep} T^i_x \cdot f'_\sigma = q_x^{-i(i-1)/2}
\on{Tr}(\Lambda^i \sigma(\on{Fr}_x)) f'_\sigma, \quad \quad i=1,\ldots,n,
\forall x\in |X|.
\end{equation} The first step in the construction of $f'_\sigma$ is
the construction of the Whittaker function.

\subsection{Whittaker functions} Introduce the following notation: for a
homomorphism $\on{Spec} R \arr X$ and an $\OO_X$--module $M$ we denote
by $M_R$ the $R$--module of sections of the pull-back of $M$ to
$\on{Spec} R$.  Denote by $\Omega$ the canonical bundle over $X$. Let
$GL_n^J(R)$ be the group of invertible $n \times n$ matrices $A =
(A_{ij})_{0\leq i,j \leq n-1}$, where $A_{ij} \in \Omega_R^{j-i}$. The
group $GL_n^J(R)$ is locally (for the Zariski topology) isomorphic to
the corresponding non-twisted group $GL_n(R)$. To establish such an
isomorphism, one has to choose a non-vanishing regular section
$\delta$ of $\Omega$, so the isomorphism is not canonical.

It is easy to see that $GL_n^J(R)$ is the group of $R$--points of a group
scheme over $X$, but in this paper we will not use this fact.

We denote by $N^J(R), T^J(R), P_1^J(R)$, etc., the corresponding
subgroups of $GL_n^J(R)$.

The twisted forms $GL_n^J(R)$ have the following advantage. Let
$u_{i,i+1}$ the $i$th component of the image of $u \in N^J(R)$ in
$N^J(R)/[N^J(R),N^J(R)]$ corresponding to the $(i,i+1)$ entry of
$u$. Then $u_{i,i+1} \in \Omega_R$.

Let us fix once and for all a non-trivial additive character $\psi: k
\arr \Ql^\times$. We define the character $\Psi_x$ of $N^J(\K_x)$ by
the formula
$$\Psi_x(u) = \prod_{i=1}^{n-1} \psi(\on{Tr}_{k_x/k}(\on{Res}_x
u_{i,i+1})),$$ and the character $\Psi$ of $N^J(\A)$ by
$$\Psi((u_x)) = \prod_{x \in |X|} \Psi_x(u_x).$$ It follows that
$\Psi(u)=1$ if $u \in N^J(\OO)$ or $u \in N^J(F)$.

For each $x\in |X|$ consider the group $GL_n^J(\K_x)$. Let $\gamma$ be
a semi-simple conjugacy class in $GL_n(\Ql)$. The following result is
due to Shintani \cite{Shi}, and Casselman and Shalika \cite{CS}.

\begin{thm}    \label{wgamma} (1) {\em There exists a unique function
$W_{\gamma,x}$ on $GL_n^J(\K_x)$ that satisfies the following properties:
\begin{itemize}
\item
$W_{\gamma,x}(gh) = W_{\gamma,x}(g), \forall h \in GL_n^J(\OO_x)$,
$W_{\gamma,x}(1)=1$;

\item
$W_{\gamma,x}(ug) = \Psi_x(u) W_{\gamma,x}(g), \forall u \in N^J(\K_x)$;

\item
$W_{\gamma,x}$ is an eigenfunction with respect to the local
Hecke-operators $T^i_x$, $i=1,\ldots,n$:
$$T^i_x \cdot
W_{\gamma,x}=q_x^{-i(i-1)/2}\on{Tr}({\Lambda}^i(\gamma))W_{\gamma,x}.$$
\end{itemize}}

(2) {\em The function $W_{\gamma,x}$ is given by the following formula.

For $\lambda=(\la_1,\ldots,\la_n) \in P^+_n$, the set of dominant weights
of $GL_n$ (i.e., such that $\la_1 \geq \la_2 \geq \ldots \geq \la_n$)
\begin{equation}    \label{cassha}
W_{\gamma,x}(\on{diag}(\pi_x^{\la_1},\ldots,\pi_x^{\la_n})) = q_x^{n(\la)}
\on{Tr}(\gamma,V(\la)),
\end{equation} where $V(\la)$ is the irreducible representation of
$GL_n(\Ql)$ of highest weight $\la$, and $n(\la) = \sum_{i=1}^n (i-1)
\la_i$.

For $\lambda \in \Z^n - P^+_n$,
$W_{\gamma,x}(\on{diag}(\pi_x^{\la_1},\ldots,\pi_x^{\la_n})) = 0$.}
\end{thm}

There is a bijection between the weight lattice of $GL_n$ and the double
quotient $N^J(\K_x)\backslash GL_n^J(\K_x)/GL_n^J(\OO_x)$, which maps
$(\la_1,\ldots,\la_n)$ to the double coset of \newline
$\on{diag}(\pi_x^{\la_1},\ldots,\pi_x^{\la_n})$. This explains the fact
that $W_{\gamma,x}$ is uniquely determined by its values at the points
$\on{diag}(\pi_x^{\la_1},\ldots,\pi_x^{\la_n})$.

\begin{rem}    \label{1} The uniqueness of the Whittaker function is
connected with the fact that an irreducible smooth representation of a
reductive group $G$ over a local non-archimedian field has at most one
Whittaker model, see \cite{GK}.

There is an explicit formula for the Whittaker function associated to an
arbitrary reductive group, due to Casselman and Shalika \cite{CS}, which we
will use in Sect.~5.\qed
\end{rem}

\medskip

Now we attach to $\sigma \in {\mathfrak A}_n$ the global Whittaker function
$W_\sigma$ on $GL_n^J(\A)$ by the formula
\begin{equation}    \label{wsigma} W_\sigma((g_x)) = \prod_{x \in |X|}
W_{\sigma(\on{Fr}_x),x}(g_x).
\end{equation} It satisfies:
\begin{itemize}
\item
$W_\sigma(gh) = W_\sigma(g), \forall h \in GL_n^J(\OO)$, $W_\sigma(1)=1$;

\item
$W_\sigma(ug) = \Psi(u) W_\sigma(g), \forall u \in N^J(\A)$; in particular,
$W_\sigma$ is left invariant with respect to $N^J(F)$;

\item For all $i=1,\ldots,n$ and $x\in |X|$ we have:
$$T^i_x \cdot W_{\sigma} = q_x^{-i(i-1)/2}\on{Tr}({\Lambda}^i
\sigma(Fr_x))W_{\sigma}.$$
\end{itemize}

\subsection{The construction of $f'_\sigma$}

Let $C^\infty(GL_n^J(\A))^{N^J(\A)}_\Psi$ be the space of
$\Ql$--valued smooth (see, e.g., \cite{BZ}) functions $f$ on
$GL_n^J(\A)$, such that $f(ug) = \Psi(u) f(g), \forall u \in N^J(\A)$;
we call such functions $(N^J(A),\Psi)$--equivariant. Let
$C^\infty(GL_n^J(\A))^{P_1^J(F)}_{\on{cusp}}$ be the space of smooth
functions $f$ on $GL_n^J(\A)$, which satisfy $f(pg) = f(g), \forall p
\in P_1^J(F)$ and are cuspidal, i.e., for each proper parabolic
subgroup of $GL_n$, whose unipotent radical $V$ is contained in $N$,
$$\int_{V^J(F) \backslash V^J(\A)} f(vg) dv = 0, \quad \quad \forall g \in
GL_n^J(\A).$$

The following result is the main step in constructing automorphic functions
for $GL_n$. The existence of the subgroup $P_1$ plays a key role in this
result, and this makes the case of $GL_n$ special.

\begin{thm}    \label{shaps}
{\em There is a canonical isomorphism
$$\phi: C^\infty(GL_n^J(\A))^{N^J(\A)}_\Psi \arr
C^\infty(GL_n^J(\A))^{P_1^J(F)}_{\on{cusp}}$$ given by the formula
$$(\phi(f))(g) = \sum_{y \in N_{n-1}^J(F)\backslash GL_{n-1}^J(F)} f
\left( \begin{pmatrix} y & 0 \\ 0 & 1 \end{pmatrix} g \right).$$ This
isomorphism commutes with the right action of $GL_n^J(\A)$ on both
spaces.}
\end{thm}

For the proof, see \cite{Sha} and \cite{PS}. Note that \thmref{shaps}
is not stated exactly in this form in \cite{Sha} but its proof can be
extracted from the proof of Theorem 5.9 there. We remark that for each
$g \in GL_n^J(\A)$ the sum above has finitely many non-zero terms. By
construction, $W_\sigma \in C^\infty(GL_n^J(\A))^{N^J(\A)}_\Psi$.  Let
$f'_\sigma = \phi(W_\sigma)$. The isomorphism of \thmref{shaps}
clearly preserves the spaces of right $GL_n^J(\OO)$--invariant
functions and commutes with the action of the Hecke operators on
them. Therefore $f'_\sigma$ is right $GL_n^J(\OO)$--invariant and
satisfies formula \eqref{hep}, i.e., it is an eigenfunction of the
Hecke operators with the same eigenvalues as those prescribed by the
Langlands conjecture. Furthermore, uniqueness of the Whittaker
function $W_\sigma$ implies that the function $f'_\sigma$ is the
unique function on $GL_n^J(\A)$ satisfying the above properties (up to
a non-zero constant multiple). Thus, the Langlands \conjref{Lan} is
equivalent to

\begin{conj}    \label{precise} For each $\sigma \in {\mathfrak G}_n$, the
function $f'_\sigma$ is left
$GL_n^J(F)$--invariant.
\end{conj}

\begin{rem} There is an approach to proving this conjecture using analytic
properties of the $L$--function of the Galois representation $\sigma$, see
\cite{JL,JPSS,La0}, which will not be discussed here.\qed
\end{rem}

\subsection{Interpretation in terms of vector bundles}    \label{geomint}

We begin by fixing notation. Recall that $B$ is the Borel subgroup in
$GL_n$ (consisting of upper triangular matrices) and $T \subset B$ is
the maximal torus (consisting of diagonal matrices). Let $P$ be the
maximal parabolic subgroup of $GL_n$ containing the subgroup $P_1$ of
$GL_n$ defined by formula \eqref{pp}.

Denote $T^J(\K_x)^+ = T^J(\K_x) \cap \on{Diag}(\OO_x)$, where
$\on{Diag}(\OO_x)$ is the set of diagonal $n \times n$ matrices with
coefficients in $\OO_x$, $B^J(\K_x)^+ = N^J(\K_x) T^J(\K_x)^+$ and
$$P^J(\K_x)^+ = \left\{ \begin{pmatrix} a & b \\ 0 & c \end{pmatrix}
\cond a \in GL_{n-1}^J(\K_x), c \in \K_x^\times, |c| \leq 1
\right\}.$$ Let $B^J(\A)^+ = \prod'_{x \in |X|} B^J(\K_x)^+$ and
$P^J(\A)^+ = \prod'_{x \in |X|} P^J(\K_x)^+$.

Denote by $Q$ the double quotient $N^J(F)\backslash B^J(\A)^+/B^J(\OO)$.
Note that since $$B^J(\A)/B^J(\OO) \simeq GL_n^J(\A)/GL_n^J(\OO),$$ $Q$ is
naturally a subset of
$N^J(F)\backslash GL_n^J(\A)/GL_n^J(\OO)$.

Let $M'_n$ denote the double quotient $P_1^J(F)\backslash
P^J(\A)^+/P^J(\OO)$. Note that since $$P^J(\A)/P^J(\OO) \simeq
GL_n^J(\A)/GL_n^J(\OO),$$ $M'_n$ is naturally a subset of
$P_1^J(F)\backslash GL_n^J(\A)/GL_n^J(\OO)$.

Finally, set $M_n = GL_n^J(F)\backslash GL_n^J(\A)/GL_n^J(\OO)$.

Let us denote by $\nu$ and $p$ the obvious projections $Q \arr M'_n$ and
$M'_n \arr M_n$, respectively.

\begin{lem}    \label{easy} {\em (1)} There is a canonical bijection
between the set $Q$ and the set of isomorphism classes of the following
data: $\{ L,(F_i),(s_i)
\}$, where $L$ is a rank $n$ vector bundle over $X$, $0=F_n \subset F_{n-1}
\subset \ldots \subset F_1 \subset F_0= L$ is a full flag of subbundles in
$L$, and $s_i: \Omega^i \arr L_i = F_i/F_{i+1}$ is a non-zero
$\OO_X$--module homomorphism.

{\em (2)} There is a canonical bijection between the set $M'_n$ and the set
of isomorphism classes of pairs $\{ L,s \}$, where $L$ is a rank $n$ vector
bundle over $X$, and $s: \Omega^{n-1} \arr L$ is a non-zero $\OO_X$--module
homomorphism.

The natural projection $\nu: Q\arr M'_n$ sends $\{ L,(F_i),(s_i) \} \in Q$
to $\{ L,s_{n-1} \} \in M'_n$.

{\em (3)} There is a canonical bijection between the set $M_n$ and the set
of isomorphism classes of rank $n$ vector bundles over $X$.

The natural map $p: M'_n \arr M_n$ corresponds to forgetting the section
$s$.
\end{lem}

\begin{proof} Recall that for a morphism $\on{Spec} R \arr X$ and an
$\OO_X$--module $M$ we denote by $M_R$ the space of sections of the
pull-back of $M$ to $\on{Spec} R$. Denote by $J^0$ the vector bundle
$\oplus_{i=0}^{n-1} \Omega^i$.

Let $Bun$ be the set of data $\{ L,\vf_{\gen},\vf_x \}$, where $L$ is
a rank $n$ bundle on $X$, and $\vf_{\gen}: J^0_F \arr L_F$ and $\vf_x:
J^0_{\OO_x} \arr L_{\OO_x}, \forall x \in |X|$ are isomorphisms
(generic and local ``trivializations'', respectively). We construct a
map $b: Bun \arr GL_n^J(\A)$ as follows. After the identifications of
$J^0_F \otimes_F \K_x \simeq J^0_{\OO_x} \otimes_{\OO_x} \K_x$ with
$J^0_{\K_x}$, and of $L_F \otimes_F \K_x \simeq L_{\OO_x}
\otimes_{\OO_x} \K_x$ with $L_{\K_x}$, $\vf_x$ and $\vf_{\gen}$ give
rise to homomorphisms $J^0_{\K_x} \arr L_{\K_x}$ which we denote by
the same characters. Let $\mu_x = (\vf_x)^{-1} \vf_{\gen}$ be the
corresponding automorphism of $J^0_{\K_x}$.

To represent the element $\mu_x$ by an $n \times n$ matrix $g_x =
(g_{x,ij})$ of the form given in Sect.~2.2, we set $g_{x,ij}$ to be
equal to the element of $\Omega^{j-i}_{\K_x}$ corresponding to the map
$\Omega^i_{\K_x} \arr J^0_{\K_x} \stackrel{\mu_x}{\longrightarrow}
J^0_{\K_x} \arr \Omega^j_{\K_x}$. Thus, $g_x$ is the
transpose\footnote{taking the transpose has some advantages, in
particular, it agrees with the conventions adopted in \cite{La1}} of
the matrix representing the action of $\mu_x$ on $J^0_{\K_x}$. The map
$b$ sends $\{ L,\vf_{\gen},\vf_x \}$ to $(g_x)_{x \in |X|} \in
GL_n^J(\A)$. It is easy to see that this map is a bijection.

Now to prove part (1) of the lemma, let us observe that given a triple
$\{ L,(F_i),(s_i) \}$, we can choose $\vf_{\gen}$ and $\vf_x$'s in
such a way that for each $j=0,\ldots,n-1$, they map
$\oplus_{i=j}^{n-1} \Omega^i_R$ to $F_{j,R} \subset L_R$ and the
associated maps $\Omega^j_R \arr F_{j,R}/F_{j+1,R}$ coincide with
restrictions of $s_j$ to $\on{Spec} R$ (here $R=F$ or $\OO_x$). With
such a choice, $g_x \in B^J(\K_x)^+, \forall x \in |X|$, and the
arbitrariness in the choice of $\vf_{\gen}$ (resp., $\vf_x$)
corresponds to left (resp., right) multiplication of $(g_x)_{x
\in |X|}$ by elements of $N^J(F)$ (resp., $B^J(\OO_x)$).

This proves part (1) of the lemma. The proof of parts (2) and (3) is
similar.
\end{proof}

Note that the function $W_\sigma$ (resp., $f'_\sigma$), which is defined on
the set \newline $N^J(F)\backslash GL_n^J(\A)/GL_n^J(\OO)$ (resp.,
$P_1^J(F)\backslash GL_n^J(\A)/GL_n^J(\OO)$) is uniquely determined by its
restriction to the subset $Q$ (resp., $M'_n$), since it is an eigenfunction
of the Hecke operators $T^n_x$.

Now $f'_\sigma$ is a function on $P_1^J(F)\backslash
GL_n^J(\A)/GL_n^J(\OO)$. Its restriction to $M'_n$, which we also denote by
$f'_\sigma$, equals, by definition, $\nu_{!}(W_\sigma)$, where $\nu_{!}$
denotes the operation of summing up a function along the fibers of the map
$\nu$ (note that these fibers are finite). \conjref{precise} can now be
stated in the following way.

\begin{conj}    \label{precise1} The function $f'_\sigma$ is constant along
the fibers of the map $p: M'_n
\arr M_n$.
\end{conj}

In the next section we discuss a geometric version of this conjecture.

\section{Conjectural geometric construction of an automorphic sheaf}

\subsection{Definitions of stacks}
Let ${\mathcal M}_n$ be the moduli stack of rank $n$ bundles on
$X$. Recall that for an $\Fq$--scheme $S$, $\on{Hom}(S,{\mathcal
M}_n)$ is the grouppoid, whose objects are rank $n$ bundles on $X
\times S$ and morphisms are isomorphisms of such bundles. Let
${\mathcal M}'_n$ be the moduli stack of pairs $\{ L,s \}$, where $L$
is a rank $n$ bundle on $X$ and $s: \Omega^{n-1} \arr {\mathcal M}'_n$
is an embedding of $\OO_X$--modules. More precisely,
$\on{Hom}(S,{\mathcal M}'_n)$ is the grouppoid, whose objects are pairs
$\{ L_S,s_S \}$, where $L_S \in \on{Ob}\on{Hom}(S,{\mathcal M}_n)$ and
 $s_S: \Omega_X^{n-1} \boxtimes \OO_S \to L_S$ is an embedding, such that the
quotient $L_S/\on{Im} s_S$ is $S$--flat; morphisms are isomorphisms of
such pairs which make the natural diagram commutative.

The set $M_n$ (resp. $M'_n$) can be identified with the set of
${\mathbb F}_q$--points of $\M_n$ (resp.  $\M'_n$). As was explained
in the introduction, we expect that $f'_\sigma$ is the function
attached to a complex of $\ell$--adic sheaves ${\mathcal S}'\si$ on
$\M'_n$. In this section we present the construction of a candidate
for the complex ${\mathcal S}'\si$ following Laumon \cite{La2}. At the
level of $\Fq$--points, this construction is actually different from
the construction of $f'_\sigma$ given in Sect.~2.3.

The reason is the following. It is easy to define a ``naive'' stack
${\mathcal Q}$ classifying triples $\{ L,(F_i),(s_i) \}$ (as in
\lemref{easy}) with ${\mathcal Q}(\Fq)=Q$ and a morphism ${\mathcal Q} \arr
{\mathcal M}'_n$ corresponding to the map of sets $\nu: Q \arr M'_n$. But
this ${\mathcal Q}$ is a disjoint union of connected components labeled by
the $n$--tuples $(d_0,\ldots,d_{n-1})$, where $d_i$ is the degree of the
divisor of zeros of the map $s_i: \Omega^i \arr F_i/F_{i+1}$. On the other
hand, the stack ${\mathcal M}'_n$ is a disjoint union of connected
components corresponding to the degree of the divisor of zeros of $s:
\Omega^{n-1} \arr L$. Recall that under $\nu$, $s_{n-1}$ becomes $s$. This
means that the fibers of $\nu$ are disconnected. Hence one can not obtain
an irreducible sheaf on ${\mathcal M}'_n$ as the direct image of a sheaf on
${\mathcal Q}$.

In this section we replace the ``naive'' stack ${\mathcal Q}$ by a
stack $\wt{\mathcal Q}$, and the Whittaker function $W_\sigma$ by a
perverse sheaf $\f\si$ on $\wt{\mathcal Q}$. The pair ($\wt{\mathcal
Q}$,$\f\si$) was first constructed by Laumon \cite{La2}.

\subsection{The stack $\wt{\mathcal Q}$}

The algebraic stack $\wt{\mathcal Q}$ is defined as follows. For an
$\Fq$--scheme $S$, $\on{Hom}(S,\wt{\mathcal Q})$ is the grouppoid,
whose objects are quintuples $\{L_S,\s_S,J_S,(J_{i,S}),(s_{i,S})\}$,
where $L_S$ and $J_S$ are rank $n$ bundles on $X \times S$, $\s_S: J_S
\arr L_S$ is an embedding of the corresponding $\OO_{X \times
S}$--modules, such that the quotient is $S$--flat, $(J_{i,S})$ is a
full flag of subbundles
$$0=J_{n,S} \subset J_{n-1,S} \subset \ldots J_{1,S} \subset J_{0,s}
=J_S,$$ and $s_{i,S}$ is an isomorphism ${\Omega}^i_X \boxtimes \OO_S
\simeq J_{i,S}/J_{i+1,S}, i=0,\ldots,n-1$. Morphisms are isomorphisms
of the corresponding $\OO_{X \times S}$--modules making all natural
diagrams commutative.

We have a natural representable morphism of stacks
$\wt{\nu}:\wt{\mathcal Q}\to\M'_n$, which for each $\Fq$--scheme $S$
maps $\{L_S,\s_S,J_S,(J_{i,S}),(s_{i,S})\}$ to the pair $\{L_S,\s_S
\circ s_{n-1,S}\}$, where $s_{n-1,S}$ is viewed as an embedding of
${\Omega}^n_X \boxtimes \OO_S$ into $J_S$.

Let $\wt{Q} = \wt{\mathcal Q}({\mathbb F}_q)$ be the set of
$\Fq$--points of $\wt{\mathcal Q}$ (see Sect.~1.7). By definition, it
consists of quintuples $\{L,\s,J,(J_i),(s_i)\}$, where $L$ and $J$ are
rank $n$ bundles on $X$, $\s:J \arr L$ is an embedding of the
corresponding $\OO_X$--modules, $(J_i)$ is a full flag of subbundles
$$0=J_n \subset J_{n-1} \subset \ldots J_1 \subset J_0 =J,$$ and
$s_i$ is an isomorphism $s_i: {\Omega}^i \simeq J_i/J_{i+1},
i=0,\ldots,n-1$.

There is a natural map of sets $r:\wt{Q}\to Q$ defined as
follows. Given an object $\{L,\s,J,(J_i),(s_i)\}$, define $F_i$ to be
the maximal locally free submodule of $L$ of rank $n-i$, which
contains the image of $J_i \subset J$ under $\mathbf{s}$. Then $(F_i)$
is a full flag of subbundles of $L$. The composition of $s_i:
{\Omega}^i \to J_i/J_{i+1}$ with the natural map $J_i/J_{i+1}\to
F_i/F_{i+1}$ induced by $\s$ is an $\OO$--module homomorphism $s'_i:
\Omega^i \arr F_i/F_{i+1}$ for each $i=0,\ldots,n-1$. Then $\{
L,(F_i),(s'_i) \}$ is a point of $Q$. Thus we obtain a map $r: \wt{Q}
\arr Q$.

\begin{lem} \label{r1} The composition $\nu\circ r:\wt{Q}\to M'_n$
coincides with the map $\wt{\nu}$. Moreover for every function $f$ on
$\wt{Q}$ we have
$\nu_{!}(r_{!}(f))=\wt{\nu}_{!}(f)$, the integrations being taken with
respect to the canonical measures on each of the three sets.
\end{lem}

\subsection{The sheaf ${\mathcal L}_E$}

In this section we recall Laumon's construction \cite{La1} of the sheaf
${\mathcal L}_E$. Let $\Coh$ be the stack classifying torsion sheaves of
finite length on $X$, i.e., for an $\Fq$--scheme $S$, $\on{Hom}(S,\Coh)$ is
the grouppoid, whose objects are coherent sheaves ${\mathcal T}_S$ on $X
\times S$, which are finite and flat over $S$ (see
\cite{La1}\footnote{the notation used in \cite{La1} is $\Coh^0$; we
suppress the upper index $0$ to simplify notation}).

Let $\Cohn$ be the open substack of $\Coh$ that classifies torsion
sheaves that have at most $n$ indecomposable summands supported at
each point. The stack $\Cohn$ can be understood as follows. Let $K$ be
a field containing $\Fq$, and let ${\mathcal T}\in \Coh(K)$. We have a
(non-canonical) isomorphism
\begin{equation}    \label{torsion}
{\mathcal T} \simeq {\mathcal O}_{X_K}/{\mathcal O}_{X_K}(-D_1) \oplus
\dots \oplus {\mathcal O}_{X_K}/{\mathcal O}_{X_K}(-D_h),
\end{equation}
where $X_K = \on{Spec} K \times_{\on{Spec} \Fq} X$, and $D_1\geq D_2
\geq \dots \geq D_h$ is a decreasing sequence (uniquely determined by
${\mathcal T}$) of effective divisors on $X_K$. The torsion sheaf
${\mathcal T}$ belongs to $\Cohn(K)$ precisely when $h\leq n$. Let $S$
be an $\Fq$--scheme. Then a torsion sheaf ${\mathcal T}_S\in\Coh(S)$
belongs to $\Cohn(S)$ if it does so at every closed point of $S$.

The stack $\Cohn$ is a disjoint union of connected components
$\Cohn=\bigcup_{m\in{\mathbb Z}_+} \Coh_{n,m}$, where the component
$\Coh_{n,m}$ classifies torsion sheaves of degree $m$; the degree of
the torsion sheaf ${\mathcal T}$ in \eqref{torsion} is $\sum_i
\on{deg} (D_i)$.  Each $\Coh_{n,m}$ has an open substack
$\Coh_{n,m}^{rss}$ classifying regular semi-simple torsion sheaves: an
$S$--point ${\mathcal T}_S$ of $\Cohn$ is said to be a point of
$\Coh_{n,m}^{rss}$ if for any $\Fqb$--point of $S$ the corresponding
sheaf over $X\underset{\on{Spec}\Fq}\times\on{Spec}\Fqb$ is isomorphic
to $\oplus \OO_X/\OO_X(-x_i)$, where the points $x_i$ are distinct.

Let $X^{(m)}$ denote the $m$th symmetric power of $X$ and let
$X^{(m),rss}$ be the complement to the divisor of diagonals in
$X^{(m)}$. We have a smooth map $X^{(m),rss}\to
\Coh_{n,m}^{rss}$. When we make base change from $\Fq$ to $\Fqb$,
$\Coh_{n,m}^{rss}$ can be identified with the quotient of $X^{(m)}$ by
$GL_1^m$ with respect to the trivial action.

Let us consider the rank $n$ local system $E=E_\sigma$ on $X$
corresponding to the Galois representation $\sigma$. Let us also write
$\pi: X^m \arr X^{(m)}$ for the natural projection. Define the $m$th
symmetric power $E^{(m)}$ of $E$ as the sheaf of invariants of $\pi_*
E^{\boxtimes m}$ under the natural action of the symmetric group
$S_m$, i.e., $E^{(m)} = (\pi_* E^{\boxtimes m})^{S_m}$. The
restriction $E^{(m)}|_{X^{(m),rss}}$ clearly descends to a local
system ${\mathcal L}_{E,m}^0$ on $\Coh_{n,m}^{rss}$, since it does
over the algebraic closure of $\Fq$.  We define the perverse sheaf
${\mathcal L}_E$ as the sheaf on $\Coh_n$, whose restriction to each
$\Coh_{n,m}$ is the Goresky-MacPherson extension of ${\mathcal
L}_{E,m}^0$, i.e., ${\mathcal L}_E|_{\Coh_{n,m}} = j_{!*}{\mathcal
L}_{E,m}^0$, where $j: \Coh_{n,m}^{rss} \hookrightarrow\Coh_{n,m}$.

We will need an explicit description of the function $L_E$
corresponding to ${\mathcal L}_E$, on the set $Coh_{n,m}$. Let
$x\in|X|$ be a closed point with residue field $k_x$. Then we can
regard $x$ as a $k_x$--rational point of $X$. Recall from Sect.~2.1
that we denote by $q_x$ the cardinality of $k_x$, and that
$q_x=q^{\deg x}$. We denote by $\Coh_{n,m}(x)$ the algebraic stack
over $k_x$ that classifies torsion sheaves of degree $m$ on
$X\underset{\on{Spec}\Fq}\times\on{Spec} k_x$ supported at $x$ that
have at most $n$ indecomposable summands. Obviously, $\Coh_{n,m}(x)$
is a locally closed sub-stack of
$\Coh_{n,m}\underset{\on{Spec}\Fq}\times\on{Spec} k_x$; we denote by
$I_{m,x}$ the corresponding embedding. Let ${\mathcal L}_{E,m,x}$ be
the pull-back to $\Coh_{n,m}(x)$ of the sheaf ${\mathcal L}_{E,m}$
under the composition
$$\Coh_{n,m}(x)\overset{I_{m,x}}\longrightarrow
\Coh_{n,m}\underset{\on{Spec}\Fq}\times\on{Spec} k_x \to\Coh_{n,m}.$$
We will denote the corresponding function by $L_{E,m,x}$.

In what follows, we fix, once and for all, a geometric point $\bar x$
over each closed point $x \in |X|$. We denote by $E_x$ the stalk of
$E$ at $\bar{x}$.

Let $P^{++}_{n,m}$ be the set $\{ \la=(\la_1,\la_2,\ldots,\la_n) \cond
\la_i \in \Z, \la_1\geq \la_2\geq \ldots\geq \la_n \geq 0,
\sum_{i}\la_i=m \}$. We can consider $P^{++}_{n,m}$ as a subset of the
set $P^+_n$ of dominant weights of $GL_n$. For $\lambda\in
P^{++}_{n,m}$, we write $E_x(\la)$ for the representation of
$GL_n(\Ql)\cong GL(E_x)$ of highest weight $\la$.

The stack $\Coh_{n,m}(x)$ has a stratification by locally closed
substacks $\Coh_{n,m}^{\la}(x)$ indexed by $\la\in P^{++}_{n,m}$. The
stratum corresponding to $\la=(\la_1,\la_2,\ldots,\la_n)\in
P^{++}_{n,m}$ parametrizes torsion sheaves of the form ${\mathcal T}
\simeq \OO_X/\OO_X(-\la_1 x) \oplus\dots\oplus\OO_X/\OO_X(-\la_n
x)$. Let ${\mathcal B}_{\la,x}$ denote the intersection cohomology
sheaf associated to the constant sheaf on the stratum
$\Coh_{n,m}^{\la}(x)$.

Let ${\mathcal T}$ be an $\Fq$--valued point of $\Coh_{n,m}$ and let
$x\in|X|$ be a closed point with residue field $k_x$. The pull-back of
${\mathcal T}$ to a sufficiently small Zariski neighborhood of $x$ in
$X\underset{\on{Spec}\Fq}\times\on{Spec} k_x$ gives rise to an object
${\mathcal T}_x$ of $\Coh_{n,m_x}(x)$, i.e., to a $k_x$--rational
point of $\Coh_{n,m_x}(x)$, for some $m_x$.  We have: $\sum_{x\in|X|}
m_x\cdot\deg(x) = m$. Moreover, there is a bijection:
\begin{equation}    \label{bijec}
\Cohn(\Fq)=\prod_{x\in|X|} {}\!\!^{'} \; \; \Coh_n(x)(k_x).
\end{equation}

The explicit description of $L_E$ is
given in the following proposition. Note that the shifts in degrees
are due to the fact that the stack $\Coh_{n,m}(x)$ has dimension $-m$,
whereas the dimension of $\Coh_{n,m}$ is zero.

\begin{prop}[\cite{La1},(3.3.8)]    \label{descr}

{\em(1)} Let ${\mathcal T}$
be an $\Fq$--point of $\Coh_{n,m}$. Then, using the notation above, we have:
$$L_{E,m}({\mathcal T})=\prod_{x\in|X|} L_{E,m_x,x}({\mathcal T}_x).$$

{\em(2)} Furthermore,
$${\mathcal L}_{E,m,x} \simeq \oplus_{\la\in P^{++}_{n,m}} {\mathcal
B}_{\la,x}[m]\left( -n(\la) \right) \otimes E_x(\la),$$ where
$n(\la)=\sum_{i=1}^n (i-1) \la_i$.
\end{prop}

\begin{rem}    \label{3} Let $x$ be an $\Fq$-rational point of $X$.
Consider now the variety ${\mathcal N}_m\subset \mathfrak{gl}_m$
of nilpotent matrices. The stack $\Coh_{m,m}(x)$ is isomorphic
to the stack ${\mathcal N}_m/GL_m$, where $GL_m$ acts on ${\mathcal N}_m$
by conjugation. Let
$\pi:\wt{\mathcal N}_m \to {\mathcal N}_m$ denote the Springer resolution
and let ${{\mathcal S}p}_m= R\pi_*\Ql$ denote the Springer sheaf on
${\mathcal N}_m$. The sheaf ${{\mathcal S}p}_m$ has a natural action of the
symmetric group $S_m$. It is shown in \cite{La1} that
$${\mathcal L}_{E,m,x} \simeq \left( {{\mathcal S}p}_m \otimes
(E_x)^{\otimes m} \right)^{S_m}|_{\Coh_{n,m}}$$ (note that $\Coh_{n,m}$ is
an open substack of $\Coh_{m,m}$). Hence the function $L_E$ associated to
${\mathcal L}_E$ can be expressed via the Kostka-Foulkes polynomials, see
\cite{La1}.

However, it will be more convenient for us to use another
interpretation of the sheaf ${\mathcal L}_{E,m,x}$, via the affine
Grassmannian (see Sect.~4.2). This interpretation allows us to express
$L_E$ in terms of the Hecke algebra ${\mathcal
H}(GL_n(\K),GL_n(\OO))$, see Sect.~5.5. The fact that the two
interpretations agree is due to Lusztig \cite{Lu1}.\qed
\end{rem}

\subsection{The sheaf $\f\si$}

Define a morphism of stacks $\alpha: \wt{\mathcal Q}\to \Cohn$ that
sends a quintuple $\{L_S,\s_S,J_S,(J_{i,S}),(s_{i,S})\}$ to the
sheaf $L_S/\on{Im} \s_S$.

Now we define a morphism $\beta: \wt{\mathcal Q} \to \Gaf$, which at
the level of $\Fq$--points sends $\{L,\s,J,(J_i),(s_i)\}$ to the sum
of $n-1$ classes in
$$\Fq \simeq Ext^1(\Omega^i,\Omega^{i-1}) \simeq
Ext^1(J_i/J_{i+1},J_{i-1}/J_i)$$ that correspond to the successive
extensions
$$0\to J_{i-1}/J_i \to J_{i-1}/J_{i+1}\to J_i/J_{i+1}\to 0.$$ Given
two coherent sheaves $L$ and $L'$ on $X$, consider the stack ${\mathcal
E}xt^1(L',L)$, such that the objects of the grouppoid
$\on{Hom}(S,{\mathcal E}xt^1(L',L))$ are coherent sheaves $L''$ on $X
\times S$ together with a short exact sequence $$0 \arr L \boxtimes
\OO_S \arr L'' \arr L' \boxtimes \OO_S \arr 0,$$ and morphisms are
morphisms between such exact sequences inducing isomorphisms at the
ends. There is a canonical morphism of stacks ${\mathcal E}xt^1(L',L)
\arr Ext^1(L',L)$. We have for each $i=1,\ldots,n-1$, a natural
morphism $\beta_i: \wt{\mathcal Q} \arr {\mathcal
E}xt^1(\Omega^i,\Omega^{i-1})$, as above. Now $\beta$ is the
composition $$\wt{\mathcal Q} \arr \prod_{i=1}^{n-1} {\mathcal
E}xt^1(\Omega^i,\Omega^{i-1}) \arr \prod_{i=1}^{n-1}
Ext^1(\Omega^i,\Omega^{i-1}) \arr \Gaf^{n-1} \arr \Gaf.$$

Let ${\mathcal I}_\psi$ be the Artin-Schreier sheaf on $\Gaf$
corresponding to the character $\psi$.

Recall that the Galois representation $\sigma$ gives rise to a rank
$n$ local system $E$ on $X$ and  to the sheaf ${\mathcal L}_E$ on
$\Cohn$. Define the sheaf $\f\si$ on $\wt{\mathcal Q}$ as
$$\f\si:=\alpha^*({\mathcal L}_E)\otimes \beta^*({\mathcal I}_\psi).$$
Note that $\wt{\mathcal Q}$ is an open sub-stack in a vector bundle
over the product of $\Cohn$ and a smooth stack that classifies
extensions $J$ as above. Hence the map $\alpha$ is smooth, and $\f\si$
is the Goresky-MacPherson extension from its restriction to the open
substack $\alpha^{-1}(\Cohn^{rss})$.

\subsection{Geometric Langlands conjecture for $GL_n$} Recall that we have
a representable morphism of stacks $\wt{\nu}:\wt{\mathcal Q}\to\M'_n$
that associates to an object $\{L,\s,J,(J_i),(s_i)\}$, the pair
$\{L,\s \circ s_{n-1}\}$.

We define the complex of $\ell$--adic sheaves ${\mathcal S}'\si$ on
$\M'_n$ to be the direct image
$${\mathcal S}'\si:=\wt{\nu}_{!}(\f\si).$$ The following conjecture of
Laumon is a geometric version of \conjref{precise1}.

\begin{conj}[\cite{La2}]    \label{princ} Let $\sigma$ be in ${\mathfrak
G}_n$ and $E$ be the corresponding irreducible $\ell$--adic local system on
$X$. Then

\begin{itemize}

\item The restriction of the complex ${\mathcal S}'\si$ to each connected
component of ${\mathcal M}'_n$ is an irreducible perverse sheaf up to a
shift in degree.

\item
${\mathcal S}'\si\simeq p^*({\mathcal S}\si)$, where $p$ is the natural
morphism $\M'_n\to \M_n$, and ${\mathcal S}\si$ is a complex of sheaves on
$\M_n$, whose restriction to each connected component of ${\mathcal M}_n$
is an irreducible perverse sheaf up to a shift.

\item The sheaf ${\mathcal S}_E$ is  an eigensheaf of the Hecke
correspondences in the sense of \cite{La1}, (2.1.1).
\end{itemize}
\end{conj}

If this conjecture is true, then the function on $M_n$ associated to the
sheaf ${\mathcal S}_E$ is the automorphic function $f_\sigma$ corresponding
to $\sigma$. The sheaf ${\mathcal S}_E$ can therefore be called the
automorphic sheaf corresponding to $\sigma$.

The conjecture means that the sheaf ${\mathcal S}'_E$ is constant along the
fibers of the morphism $p$. Thus, it is analogous to
\conjref{precise1}. The advantage of dealing with \conjref{princ} as
compared to \conjref{precise1} is that while the latter is a global
statement, one could use local geometric information about the sheaf
${\mathcal S}'_E$ to tackle \conjref{princ} (as Drinfeld did in the case of
$GL_2$ \cite{Dr}).

\begin{rem} The above conjecture is obviously false if one does not assume
the irreducibility of $\sigma$ (the complex ${\mathcal S}'\si$ must be
corrected by the corresponding ``constant terms'' in this
case). However one can construct the automorphic sheaves ${\mathcal
S}\si$ corresponding to $\sigma$'s, which are direct sums of
one-dimensional representations, by means of the geometric Eisenstein
series \cite{La:e}.\qed
\end{rem}

\subsection{Main theorem} Let $S'\si$ denote the function on $M'_n$
associated to ${\mathcal S}'\si$. If \conjref{princ} is true, then
this function has the same properties as the function $f'_\sigma$
defined in Sect.~2.3. Recall that these properties uniquely determine
$f'_\sigma$ up to a non-zero factor. Therefore \conjref{princ} can be
true only if the functions $S'_E$ and $f'_\sigma$ are
proportional. This was conjectured by Laumon in \cite{La2} (Conjecture
3.2). One of our motivations was to prove this conjecture. More
precisely, we prove the following:

\begin{thm}    \label{prin} {\em The functions $S'\si$ and $f'_\sigma$ are
equal.}
\end{thm}

\thmref{prin} means that the function $f'_\sigma$ does come from a
complex of $\ell$--adic sheaves on ${\mathcal M}'_n$. It also provides
a consistency check for \conjref{princ}. Laumon has proved
\thmref{prin} in \cite{La1} for $GL_2$ by a method different from the
one we use below. We derive \thmref{prin} for $GL_n$ with arbitrary
$n$ from the following statement.

\begin{prop} \label{r} Let $F\si$ be the function on $\wt{Q}$ corresponding
to $\f\si$. The function $r_{!}(F\si)$ coincides with the restriction
of the Whittaker function $W_\sigma$ to $Q$.
\end{prop}

\thmref{prin} immediately follows from \propref{r} because of \lemref{r1}
as shown on the diagram below.

\setlength{\unitlength}{1mm}

\begin{center}
\begin{picture}(65,60)(-30,-55)
\put(-4,0){\vector(-1,0){20}}
\put(4,0){\vector(1,0){20}}
\put(0,-3.5){\vector(0,-1){20}}
\put(0,-32){\vector(0,-1){20}}
\put(2.5,-2.5){\vector(1,-1){8.5}}
\put(11,-15.5){\vector(-1,-1){8.5}}
\put(-1.5,-1){$\widetilde{Q}$}
\put(-30,-1){${\mathbb F}_q$}
\put(-2,-29){$M'_n$}
\put(-2,-56){$M_n$}
\put(12,-14){$Q$}
\put(26,-1){${Coh}_n$}
\put(-15,2){$\beta$}
\put(13,2){$\alpha$}
\put(-4,-14){$\widetilde{\nu}$}
\put(-4,-43){$p$}
\put(8,-6){$r$}
\put(8,-21){$\nu$}
\end{picture}
\end{center}

\bigskip

The proof of \propref{r} will occupy Sects.~4 and 5 below.

\propref{r} can be interpreted in the following way. Let $Coh_n$
(resp., $Coh_{n,m}(x)$) be the set of $\Fq$--points (resp.,
$k_x$--points) of $\Coh_n$ (resp., $\Coh_{n,m}(x)$). Set $Coh_n(x) =
\cup_{m\geq 0} Coh_{n,m}(x)$. Then by formula \eqref{bijec}, $Coh_n =
\prod'_{x\in |X|} Coh_n(x)$. We have:
$$Coh_n(x) \simeq P^{++}_n = \cup_{m\geq 0} P^{++}_{n,m}.$$ Hence we
can identify $Coh_n(x)$ with the set $$\{
\on{diag}(\pi_x^{\la_1},\ldots,\pi_x^{\la_n})|\la_1 \geq \ldots \geq
\la_n \geq 0 \}.$$ The Whittaker function $W_\sigma$ can then be
restricted to the $Coh_n$, and it is uniquely determined by this
restriction, see Sect.~2.2. Thus, both $L_E$ and $W_\sigma$ give rise
to functions on $Coh_n$. For each point $t \in Coh_n$, the value of
$L_E$ at $t$ is given by taking the alternating sum of traces of the
Frobenius on the stalk cohomologies of ${\mathcal L}_E$ at $t$. On the
other hand, the value of $W_\sigma$ is given by the trace of the
Frobenius on the top stalk cohomology of ${\mathcal L}_E$ at $t$ (see
\cite{La1}). Therefore \propref{r} says that the contributions of all
stalk cohomologies, other than the top one, are killed by the
summation along the fibers of the projection $r$ against the
non-trivial character $\Psi$.

\begin{rem} As we mentioned in the introduction, Drinfeld has proved a
version of
\conjref{princ} for $GL_2$. The case of $GL_2$ is special as explained
below.

Let $\wt{\mathcal M}'_2$ be the open substack of ${\mathcal M}'_2$, which
parametrizes $\{ L,s \}$, such that the image of $s: \Omega \arr L$ is a
maximal invertible subsheaf of $L$. Due to the Hecke eigenfunction property
of $f'_\sigma$ with respect to $T^2_x$, $f'_\sigma$ is uniquely determined
by its restriction to $\wt{M}'_2 = \wt{\mathcal M}'_2(\Fq)$.

But the map $r$ is a bijection over $\nu^{-1}\left( \wt{M}'_2 \right)$, and
hence $\nu^{-1} \left( \wt{M}'_2 \right)$ can be considered as a subset of
$\wt{Q}$. Clearly, the map $\alpha$ of Sect.~3.3 sends $\nu^{-1} \left(
\wt{M}'_2 \right)$ to $Coh_1 \subset Coh_2$. Furthermore, the restriction
of ${\mathcal L}_E$ to the stack $\Coh_1$ is simply a sheaf, i.e., it has
stalk cohomology in only one degree. Therefore on $Coh_1$ the function
$L_E$ equals the Whittaker function $W_\sigma$. Hence, restricted to
$\wt{\mathcal M}'_2$, the geometric construction coincides with the
construction described in Sect.~2.3.\qed
\end{rem}

\section{Reduction to a local statement}
\subsection{Adelic interpretation of $\wt{Q}$}

Recall that $\wt{Q}$ is the set of isomorphism classes of ${\mathbb
F}_q$--points of $\wt{\mathcal Q}$, and in Sect.~3.1 we defined a map
$r:\wt{Q} \to Q$. Recall further that $Coh_n(x) \simeq P^{++}_n
=\cup_{m\geq 0} P^{++}_{n,m}$, and $Coh_n = \prod'_{x\in |X|}
Coh_n(x)$.

Denote $GL_n^J(\K_x)^+ = GL_n^J(\K_x) \cap \on{Mat}_n^J(\OO_x)$. Let
$GL_n^J(\A)^+$ be the restricted product $\prod'_{x\in |X|}
GL_n^J(\K_x)^+$. We have bijections:
\begin{equation}    \label{iso1} GL_n^J(\OO_x) \backslash
GL_n^J(\K_x)^+/GL_n^J(\OO_x) \simeq Coh_n(x).
\end{equation} and
\begin{equation}    \label{iso2} GL_n^J(\OO) \backslash
GL_n^J(\A)^+/GL_n^J(\OO) \simeq Coh_n.
\end{equation}

\begin{prop}    \label{adelic} {\em (1)} There is a bijection between the
set $\wt{Q}$ and the set
$$(N^J(F)\backslash N^J(\A))\underset{N^J(\OO)}\times
(GL_n^J(\A)^+/GL_n^J(\OO)).$$ The group $N^J(\OO)$ acts on the product
$(N^J(F)\backslash N^J(\A))\times (GL_n^J(\A)^+/GL_n^J(\OO))$ according to
the rule $y\cdot (u,g)=(u\cdot y^{-1},y\cdot g)$.

{\em (2)} The map $r: \wt{Q} \arr Q$ identifies with the map
$$(N^J(F)\backslash N^J(\A))\underset{N^J(\OO)}\times
(GL_n^J(\A)^+/GL_n^J(\OO)) \to Q \subset N^J(F)\backslash
GL_n^J(\A)/GL_n^J(\OO)$$ given by $(u,g)\to u\cdot g$.

The map $\al$ sends $(u,g)$ to the image of $g$ in $GL_n^J(\OO) \backslash
GL_n^J(\A)^+/GL_n^J(\OO) \simeq Coh_n$.

The map $\beta: \wt{Q} \arr \Fq$ is the composition of the natural map
$N^J(F) \backslash N^J(\A)/N^J(\OO) \arr \Fq^{n-1}$ and the summation
$\Fq^{n-1} \arr \Fq$.
\end{prop}

\begin{proof}
We will use the notation introduced in the proof of \lemref{easy}. Let
$\{ L,\s,J,(J_i),$ $(s_i) \}$ be an element of $\wt{Q}$. Then the triple
$\{ J,(J_i),(s_i) \}$ is an element of $Q$. Hence we can associate to
it homomorphisms $\vf^J_x: J^0_{\K_x} \arr J_{\K_x}$, $\vf^J_{\gen}:
J^0_{\K_x} \arr J_{\K_x}$ and $\mu^J_x = (\vf^J_x)^{-1} \vf^J_{\gen}:
J^0_{\K_x} \arr J^0_{\K_x}$, as in the proof of \lemref{easy}.

On the other hand, let $\vf_x^L$ be an isomorphism $J^0_{\OO_x} \arr
L_{\OO_x}$. We extend it to a homomorphism $J^0_{\K_x} \arr L_{\K_x}$,
which we denote by the same symbol. Denote by $\s_x$ the homomorphism
$J_{\K_x} \arr L_{\K_x}$ induced by $\s$. Consider the automorphism $\nu_x
= (\vf_x^L)^{-1} \s_x \vf_x^J$ of $J^0_{\K_x}$.

Now we assign to $\{ L,\s,J,(J_i),(s_i) \}$ the element
$((u_x),(g_x))$ of $GL_n^J(\A) \times GL_n^J(\A)$, where $u_x$ is the
transpose of the matrix representing $\mu_x$, and $g_x$ is the
transpose of the matrix representing $\nu_x$ (see the proof of
\lemref{easy}). By construction, $u_x \in N^J(\K_x)$ and $g_x \in
GL_n^J(\K_x)^+$. Furthermore, the arbitrariness in the choice of
$\vf_{\gen}^J$ corresponds to left multiplication of $u_x$ by elements
of $N^J(F)$, the arbitrariness in $\vf_x^L$ corresponds to right
multiplication of $g_x$ by elements of $GL_n^J(\OO_x)$, and the
arbitrariness in $\vf_x^J$ corresponds to the action of $N^J(\OO_x)$
on $(u_x,g_x)$ according to the rule $y \cdot (u_x,g_x) = (u_x \cdot
y^{-1}, y \cdot g_x)$.

This proves part (1) of the proposition. The proof of part (2) is now
straightforward.
\end{proof}

Recall that $L_E$ is the function on $Coh_n$ associated to the sheaf
${\mathcal L}_E$. Denote by $L_{E,x}$ the function on $Coh_n(x)$,
whose restriction to $\Coh_{n,m}(x)(\Fq)$ is the function associated
to ${\mathcal L}_{E,m,x}$. Part (1) of \propref{descr} implies that
$$L_E((\la_x)) = \prod_{x \in |X|} L_{E,x}(\la_x),$$ for all $(\la_x)
\in \prod'_{x \in |X|} P_n^{++}$. Using the bijection \eqref{iso1}
(resp., \eqref{iso2}) we consider $L_{E,x}$ (resp., $L_E$) as a
function on $GL_n^J(\K_x)^+$ (resp., $GL_n^J(\A)^+$). Let $\Psi_x:
N^J(\K_x)\to\Ql^{\times}$ be the character defined in Sect.~2.2. Note
that the function $\Psi_x$ (resp., $L_{E,x}$) is right (resp., left)
$N^J(\OO_x)$--invariant.

Recall further that $F_E$ is the function on $\wt{Q}$ associated to
the sheaf ${\mathcal F}_E$. We conclude:

\begin{lem}    \label{ide} Under the isomorphism of \propref{adelic},
$$F_E = \prod_{x\in |X|} \Psi_x \times L_{E,x}.$$
\end{lem}

Let us extend the function $L_{E,x}$ by zero from $GL_n^J(\K_x)^+$ to
$GL_n^J(\K_x)$. Then \propref{adelic} and \lemref{ide} imply:
$$(r_! F_E)((g_x)) = \prod_{x \in |X|} \; \; \sum_{u_x \in
N^J(\K_x)/N^J(\OO_x)} L_{E,x}(u_x^{-1} \cdot g_x) \Psi_x(u_x),$$ for
all $g_x \in GL_n^J(\K_x)^+$ (each sum is actually finite). Let $du_x$
be the Haar measure on $N^J(\K_x)$ normalized so that
$\int_{N^J(\OO_x)} du_x = 1$. Using left $GL_n^J(\OO_x)$--invariance
of $L_{E,x}$, we can rewrite the last formula as
\begin{equation}    \label{globalint} (r_! F_E)((g_x)) = \prod_{x \in |X|}
\int_{N^J(\K_x)} L_{E,x}(u_x^{-1} \cdot g_x) \Psi_x(u_x) \; du_x.
\end{equation}

\propref{r} states that $(r_! F_E)((g_x)) =
W_\sigma((g_x))|_Q$. According to formulas \eqref{wsigma} and
\eqref{globalint}, this is equivalent to the formula
\begin{equation}    \label{show}
\int_{N^J(\K_x)} L_{E,x}(u_x^{-1} \cdot g_x) \Psi_x(u_x) \; du_x =
W_{\sigma(\on{Fr}_x)}(g_x),
\end{equation}
for all $g_x \in GL_n^J(\K_x)^+$. Since both the left and the right
hand sides of \eqref{show} are left $(N^J(\K_x),\Psi_x)$--equivariant
and right $GL_n^J(\OO)$--invariant, it suffices to check formula
\eqref{show} when $g_x = \on{diag}(\pi^{\nu_1},\ldots,\pi^{\nu_n})$,
where $\nu = (\nu_1,\ldots,\nu_n) \in P^{++}_n$.

Using the explicit formula \eqref{cassha}, we reduce \propref{r} to the
following local statement.

\begin{prop}    \label{Fplus}
\begin{equation}    \label{oh}
\int_{N^J(\K_x)} L_{E,x}(u_x \cdot
\on{diag}(\pi_x^{\nu_1},\ldots,\pi_x^{\nu_n})) \Psi^{-1}(u_x) \; du_x
= q_x^{n(\nu)} \on{Tr}(\sigma(\on{Fr}_x),E_x(\nu)),
\end{equation} where $\nu \in P^{++}_n$.
\end{prop}

\subsection{Positive part of the affine Grassmannian}

Recall that $J^0 = \oplus_{i=0}^{n-1} \Omega^i$.  From now on we work
in the local setting. Hence we choose once and for all a
trivialization of $J^0$ on the formal neighborhood of $x \in |X|$ and
identify $GL_n^J(\K_x)$ with $GL_n(\K_x)$. For this reason we suppress
the index $J$ in what follows.

According to part (2) of \propref{descr},
\begin{equation}    \label{le} L_{E,x} = \sum_{\la \in P^{++}_n}
\on{Tr}(\sigma(\on{Fr}_x),E_x(\la)) \cdot B_\la,
\end{equation} where $B_\la$ is the function on $Coh_n(x)$ associated
to the sheaf ${\mathcal B}_{\la,x}[m](-n(\la))$ (see Sect.~3.2). We
view $B_\la$ as a $GL_n(\OO_x)$--invariant function on
$GL_n(\K_x)/GL_n(\OO_x)$.

Let $x$ be a closed point of $X$. To simplify notation, from now on we
will assume that $x$ is an $\Fq$--rational point of $X$ (otherwise, we
simply make a base change from $\Fq$ to the residue field $k_x$ of
$x$).

Define the functor which sends an $\Fq$--scheme $S$ to the set of
isomorphism classes of pairs $\{ L,t \}$, where $L$ is a rank $n$
bundle on $X \times S$ and $t$ is its trivialization on $(X \times S)
- (\{ x \} \times S)$. This functor is representable by an ind--scheme
$\G_x$ (see \cite{BL1}), which we call the affine Grassmannian $\G_x$
(for the group $GL_n$). The ind--scheme $\G_x$ splits into a disjoint
union of connected components: $\G_x=\cup_{m\in{\mathbb Z}} \G_x^m$
indexed by the degree of $L$ for $\{ L,t \} \in \G_x$.

There is a bijection between the set $Gr_x$ of ${\mathbb F}_q$--points
of $\G_x$ and the quotient $GL_n(\K_x)/GL_n(\OO_x)$, see
\cite{BL}. The analogous quotient over the field of complex numbers is
known as the affine, or periodic, Grassmannian. This explains the name
that we use.

There exists a proalgebraic group $\GG(\OO_x)$ whose set of
$\Fq$--points is $GL_n(\OO_x)$. The group $\GG(\OO_x)$ acts on $\G_x$,
and its orbits stratify $\G_x$ by locally closed finite-dimensional
subschemes $\G_x^\la$ indexed by the set $P^+_n$ of dominant weights
$\la$ of $GL_n$. The stratum $\G_x^\la$ is the $\GG(\OO_x)$--orbit of
the coset $\on{diag}(\pi_x^{\la_1},\ldots,\pi_x^{\la_n}) \cdot
GL_n(\OO_x) \in Gr_x$.

Recall that there is an inner product on the set of $GL_n$ weights
defined by the formula $(\la,\mu) = \sum_{i=1}^n \la_i \mu_i$, and
that $\dim \G_x^\la = 2(\la,\rho)$, where $\rho$ is the half sum of
the positive roots of $GL_n$, $\rho =
((n-1)/2,(n-3)/2,\ldots,-(n-1)/2)$. Let $\ol{\mathbb Q}_{\ell,\la}$ be
the constant sheaf supported on the stratum $\G_x^\la$. Denote by
${\mathcal A}_{\la,x}$ the intersection cohomology sheaf on the
closure of $\G_x^\la$, which is the Goresky-MacPherson extension of
the sheaf $\ol{\mathbb Q}_{\ell,\la}[2(\la,\rho)]((\la,\rho))$. Let
$A_\la$ be the $GL_n(\OO_x)$--invariant function on
$GL_n(\K_x)/GL_n(\OO_x)$, which is the extension by zero of the
function associated to ${\mathcal A}_{\la,x}$.

Define now the closed subscheme $\G_x^+$ of $\G_x$ which at the level
of points corresponds to pairs $\{L,t\}$, for which $t$ extends to an
embedding of $\OO_X$--modules $\OO_X^{\oplus n} \arr L$. The scheme
$\G_x^+$ also splits into a disjoint union of connected components
$\G_x^+=\cup_{m\in \Z^+} \G_x^{m,+}$, where $\G_x^{m,+} = \G_x^+ \cap
\G_x^m$. It is clear that $\G_x^{m,+}$ is a union of the strata
$\G_x^\la$, for $\la\in P^{++}_{n,m}$ (see Sect.~3.2). The set
$\G_x^+(\Fq)$ identifies with the quotient $GL_n(\K_x)^+/GL_n(\OO_x)$.

Consider the morphism $q_{m,x}: \G_x^{m,+} \to \Coh_{m,n}(x)$, 
which sends a pair $\{L,t\}$ to the quotient  $\OO_X^{\oplus
n}/\on{Im} t^*$, where $t^*: L^* \arr \OO_X^{\oplus n}$. This morphism
can be described as follows. There is a natural vector bundle
$\wt{\Coh}_{m,n}$ over the stack $\Coh_{m,n}$, whose fiber at
${\mathcal T}$ is $\on{Hom}(\OO_X^{\oplus n},{\mathcal T})$;
$\G_x^{m,+}$ is an open substack of the total space of the bundle
$\wt{\Coh}_n$, which corresponds to epimorphic elements of
$\on{Hom}(\OO_X^{\oplus n},{\mathcal T})$. The map $q_{m,x}$ is simply
the projection from this substack to the base. Hence $q_{m,x}$ is a
smooth morphism of algebraic stacks. It is clear that it preserves the
stratification (compare with \cite{Lu1}). Note that $\dim \G_x^{m,+} =
m(n-1)$ and $\dim \Coh_{n,m}(x) = -m$. Therefore

\begin{lem} \label{p}
\begin{equation}    \label{smooth} q_{m,x}^* {\mathcal
B}_{\la,x}[m](-n(\la)) = {\mathcal A}_{\la,x}[-m(n-1)](-(\la,\rho)-n(\la)).
\end{equation}
\end{lem}

Recall that $m=|\la|=\sum_{i=1}^n \la_i$. The lemma implies that as
functions on \newline $GL_n(\K_x)^+/GL_n(\OO_x)$,
\begin{equation}    \label{equ} B_\la = (-1)^{|\la|(n-1)}
q_x^{(\la,\rho)+n(\la)} A_\la.
\end{equation}

Hence, according to \eqref{le},
\begin{equation}    \label{le1} L_{E,x} = \sum_{\la \in P^{++}_n}
(-1)^{|\la|(n-1)} q_x^{(\la,\rho)+n(\la)}
\on{Tr}(\sigma(\on{Fr}_x),E_x(\la)) \cdot A_\la,
\end{equation}

This formula will be used in the next section in the proof of
\propref{Fplus}.

\section{Proof of the local statement} In this section we will state and
prove a general result concerning reductive groups over local
non-archimedian field of positive characteristic. In the case of $GL_n$
this result implies
\propref{Fplus}.

\subsection{General set-up} Let $G=G({\mathcal K})$ be a connected,
reductive, split algebraic group over the field ${\mathcal
K}=\Fq((\pi))$, $T$ -- its split maximal torus contained in a Borel
subgroup $B$, and $N$ -- the unipotent radical of $B$. We again denote
by ${\mathcal O}$ the ring of integers of ${\mathcal K}$, by $\pi$ a
generator of its maximal ideal, and by $q$ the cardinality of the
residue field $k={\mathbb F}_q$. Let $K$ be the compact subgroup
$G({\mathcal O})$ of $G$. We fix a Haar measure of $G$, such that $K$
has measure $1$.

Let ${\mathcal H}(G,K)$ denote the Hecke algebra of $G$ with respect
to $K$, i.e., ${\mathcal H}(G,K)$ is the algebra of $\Ql$--valued
compactly supported $K$--bi-invariant functions on $G$ with the
convolution product:
\begin{equation}    \label{conv} (f_1 \cdot f_2)(g) = \int_G f_1(x)
f_2(gx^{-1}) \; dx.
\end{equation}

Let $\GL$ be the Langlands dual group of $G$ (without the Weil group
of $\K$), and $\PL$ (resp., $\PL^+$) be the set of weights (resp.,
dominant weights) of $\GL$. Each $\la$ can be viewed as a
one-parameter subgroup of $G(\K)$, and hence $\la(\pi)$ is a
well-defined element of $G(\OO_x)$. We denote by $c_\la$ the
characteristic function of the double coset $K \la(\pi) K \subset
G$. The functions $c_\la$ form a basis of ${\mathcal H}(G,K)$.

Let $\on{Rep} \GL(\Ql)$ denote the Grothendieck ring of the category
of finite-dimensional representations of $\GL(\Ql)$. We consider it as
a $\Ql$--algebra. If $V$ is a finite-dimensional representation of
$\GL$, denote by $[V]$ the corresponding element of $\on{Rep}
\GL(\Ql)$. In particular, for each $\la \in \PL^+$, let $V(\la)$ be
the finite-dimensional representation of $\GL$ with highest weight
$\la$.

The following statement, often referred to as the Satake isomorphism, is
well-known, see \cite{Sa,La,Mac,G}.

\begin{thm}    \label{Satake} {\em There is a unique isomorphism $S:
\on{Rep} \GL(\Ql) \arr {\mathcal H}(G,K)$, which maps $[V(\la)]$ to}
$$H_\la =  q^{-(\la,\rho)} \left( c_\la + \sum_{\mu \in \PL^+;
\mu<\la} a_{\la\mu} c_\mu \right), \quad \quad a_{\la\mu} \in \Z.$$
\end{thm}

\begin{rem}    \label{character} Each semi-simple conjugacy class $\gamma$
of the group $\GL(\Ql)$ defines a homomorphism $\chi_\ga: \on{Rep} \GL(\Ql)
\arr \Ql$, which maps $[V]$ to $\on{Tr}(\gamma,V)$. We denote the
corresponding homomorphism ${\mathcal H} \arr \Ql$ by the same symbol
$\chi_\ga$. This allows us to identify the spectrum of the commutative
algebra ${\mathcal H}$ with the set of semi-simple conjugacy classes of
$\GL(\Ql)$. In particular, we have: $\chi_\ga(H_\la) =
\on{Tr}(\ga,V(\la))$.\qed
\end{rem}

\subsection{Hecke algebra and the affine Grassmannian}    \label{he}

Let again $X$ be as in Sect.~1.1, and $x$ be its $\Fq$--point. We
define, in the same way as in Sect.~4.2 for $G=GL_n$, the ind--scheme
${{\mathcal G}r}(G)=\G(G)_x$ that classifies pairs $\{{\mathcal
P},t\}$, where ${\mathcal P}$ is a principal $G$--bundle on $X$ and
$t$ is its trivialization over $X-x$. The ind--scheme structure on
$\G(G)$ is described, e.g., in \cite{LS}. Note that due to the results
of \cite{BL,DS}, the global curve $X$ is inessential in the above
definition. We could simply take $X = \on{Spec} \OO]$. In particular,
there is a bijection between the set of ${\mathbb F}_q$--points of
$\G(G)$ and the set $G/K$.

There is a proalgebraic group $\GG(\OO)$, whose set of $\Fq$--points
is $G(\OO)$. This group acts on $\G(G)$, and its orbits stratify
$\G(G)$ by locally closed finite-dimensional subschemes $\G(G)^\la$
indexed by the set $\PL^+$ of dominant weights $\la$ of $\GL$. The
stratum $\G(G)^\la$ is the $G(\OO)$--orbit of the coset $\la(\pi)
\cdot G(\OO)$, where $\la(\pi) \in T(\K) \subset G(\K)$ is defined
above.

Denote by $(\la,\rho)$ the pairing between $\la \in \PL$ and the sum
of the fundamental coweights $\rho$ of $\GL$. Let $\ol{\mathbb
Q}_{\ell,\la}$ be the constant sheaf supported on the stratum
$\G(G)^\la$. Denote by ${\mathcal A}_\la={\mathcal A}_{\la,x}$ the
intersection cohomology sheaf on the closure of $\G(G)^\la$, which is
the Goresky-MacPherson extension of the sheaf $\ol{\mathbb
Q}_{\ell,\la}[2(\la,\rho)]((\la,\rho))$ (note that $\dim \G(G)^\la =$
$2(\la,\rho)$). Let $A_\la$ be the function associated to ${\mathcal
A}^\la$. We use the same notation for its extension by zero to the
whole $\G(G)$. Clearly, the functions $A_\la, \la \in \PL^+$, form a
basis in the $\Ql$--vector space $\Ql(G/K)^K$ of $K$--invariant
functions on $G/K$ with compact support. We have an isomorphism of
vector spaces: ${\mathcal H} \simeq \Ql(G/K)^K$, which commutes with
the action of ${\mathcal H}$. Therefore $H_\la \in {\mathcal H}, \la
\in \PL$, can also be considered as elements of $\Ql(G/K)^K$.

\begin{prop}    \label{hla}
$H_\la = (-1)^{2(\la,\rho)} A_\la$.
\end{prop}

This result is due to Lusztig \cite{Lu1,Lu2} and Kato \cite{Ka} (see, e.g.,
Theorem 1.8, Lemma 2.7, formula (3.5) of \cite{Ka}). It implies that
\begin{equation}    \label{stalk} H_\la(y) = (-1)^{2(\la,\rho)} \sum_{i \in
\Z} \dim \on{H}^i({\mathcal A}_\la)|_y \; q^{i/2},
\end{equation} where $\on{H}^j({\mathcal A}_\la)|_y$ is the $j$th
stalk cohomology of
${\mathcal A}_\la$ at $y \in Gr(G)$.

\subsection{Fourier transform} Let us denote by $\on{Res}: \K \arr \Fq$ the
map defined by the formula
$$\on{Res}\left( \sum_{n \in \Z} f_i \pi^i \right) = f_{-1}.$$  We define a
character $\Psi$ of $N$ in the following way:
$$\Psi(u) = \sum_{i=1}^l \psi\left(\on{Res}(u_i)\right),$$ where $u_i,
i=1\ldots,l=\dim N/[N,N]$ are natural coordinates on $N/[N,N]$
corresponding to the simple roots and $\psi: \Fq \arr \Ql^{\times}$ is
a fixed non-trivial character.

Consider the space $\Ql(G/K)^N_\Psi$ of left $(N,\Psi)$--equivariant and
right $K$--invariant functions on $G$ that have a compact support modulo
$N$. This space is a module over ${\mathcal H}(G,K)$, with the action
defined by formula \eqref{conv}. For $\la \in \PL^+$, let $\phi_\la$ be the
function from $\Ql(G/K)^N_\Psi$, which vanishes outside the $N$--orbit of
$\la(\pi)$ and equals $q^{-(\la,\rho)}$ at $\la(\pi)$. The elements
$\phi_\la, \la \in \PL^+$, provide a $\Ql$--basis for $\Ql(G/K)^N_\Psi$.

Define the linear map $\Phi: \Ql(G/K)^K \arr \Ql(G/K)^N_\Psi$ by
the formula
\begin{equation}    \label{check1} (\Phi(f))(g) = \int_{N} f(ug)
\Psi^{-1}(u) \; du,
\end{equation} where $du$ stands for the Haar measure on $N$ normalized so
that $\int_{N({\mathcal O})} \; du = 1$.

\begin{lem}    \label{F} The map $\Phi$ defines an isomorphism of
${\mathcal H}$--modules $$\Ql(G/K)^K \simeq \Ql(G/K)^N_\Psi.$$
\end{lem}

\begin{proof} Let $c_\la \in \Ql(G/K)^K$ be the characteristic function of
the $K$--orbit of the coset $\la(\pi) \cdot K \in G/K$. It follows from the
definition that $$\Phi(c_\la) = q^{(\la,\rho)} \phi_\la + \sum_{\mu \in
\PL^+; \mu<\la} b_{\la,\mu} \phi_\mu.$$ This implies the lemma.
\end{proof}

It is natural to call $\Phi$ the Fourier transform. We are now ready to
state our main local theorem.

\begin{thm}    \label{local} {\em The map $\Phi$ sends $H_\la$ to
$\phi_\la$.}
\end{thm}

This theorem is equivalent to the formula
\begin{equation}    \label{formula}
\int_{N} H_\la(u\cdot \nu(\pi)) \Psi^{-1}(u) \; du = q^{-(\la,\rho)}
\delta_{\la,\nu}.
\end{equation}

\subsection{Proof of \thmref{local}}    \label{proof} Our proof relies on
the result of Casselman-Shalika (and Shintani for $G=GL_n$), which
describes the values of the Whittaker function at the points
$\mu(\pi)$ (cf. Theorem 2.1 and \remref{1}).

Let $\gamma$ be a semi-simple conjugacy class in $\GL$ and let
$W_{\gamma}$ be a $\Ql$-valued function on $G$ with the following three
properties:

\begin{itemize}
\item
$W_{\gamma}(gh) = W_{\gamma}(g), \forall h \in K$,
$W_{\gamma}(1)=1$;

\item
$W_{\gamma}(ug) = \Psi^{-1}(u) W_{\gamma}(g), \forall u \in N$;

\item
\begin{equation}    \label{hecke}
\int_G f(x) W_{\gamma}(gx) dx = \chi_\gamma(f) W_{\gamma}(g), \quad \quad
\forall g \in G, f \in {\mathcal H}
\end{equation} (see \remref{character} for the definition of $\chi_\gamma$).
\end{itemize}

\begin{thm}[\cite{CS},\cite{Shi}]    \label{prime} {\em The function
$W_\gamma$ satisfying these properties exists, and it is unique. For $\mu
\in \PL$, the value of this function at $\mu(\pi)$ is
\begin{equation}    \label{css} W_{\gamma}(\mu(\pi))=q^{-(\mu,\rho)}
\on{Tr}(\gamma,V(\mu)),
\end{equation} if $\mu$ is a dominant weight, and $0$, otherwise.}
\end{thm}

The function $W_\ga$ is called the Whittaker function corresponding to
$\ga$.

Now we can prove \thmref{local}. Let $\gamma$ be as above and let
$s_{\gamma}$ be a linear functional $\Ql(G/K)^N_\Psi \to \Ql$ given by the
formula
\begin{equation}    \label{sgamma} s_{\gamma}(\phi)=\int_{N\backslash G}
W_{\gamma}(g) \phi(g) \; dg,
\end{equation} where $dg$ is the measure on $N \backslash G$ induced by the
Haar measure on $G$ from Sect.~5.1 and the Haar measure on $N$ from
Sect.~5.3. By construction, the function $W_{\gamma}(u) \phi(u)$ is left
$N$--invariant. The integral \eqref{sgamma} converges, because, by
definition, $\phi$ has compact support modulo $N$.

\begin{lem}   \label{mapp} The map $s_{\gamma}$ is a homomorphism of
${\mathcal H}(G,K)$--modules
$\Ql(G/K)^N_\Psi \to \ol{\mathbb Q}_{\ell,\ga}$, where $\ol{\mathbb
Q}_{\ell,\ga}$ is the one-dimensional representation of ${\mathcal H}(G,K)$
corresponding to its character $\chi_\ga$.
\end{lem}

\begin{proof} Each $f \in {\mathcal H}$ acts on $\Ql(G/K)^N_\Psi$ by
mapping $\phi \in
\Ql(G/K)^N_\Psi$ to $f \cdot \phi$. By definition of the convolution product
(see formula \eqref{conv}), we have:
$$(f \cdot \phi)(y) = \int_G f(x) \phi(yx^{-1}) \; dx.$$ Hence
$$s_\gamma(f \cdot \phi) = \int_{N\backslash G} W_{\gamma}(g) \left( \int_G
f(x) \phi(g x^{-1}) \; dx \right) \; dg.$$ Changing the order of
integration and using the invariance of the Haar measure, we obtain
$$s_\gamma(f \cdot \phi) = \int_{N\backslash G} \left( \int_G W_\gamma(gx)
f(x) \; dx \right) \phi(g) \; dg.$$ By \eqref{hecke},
$$s_\gamma(f \cdot \phi) = \chi_{\gamma}(f) s_\gamma(\phi).$$
\end{proof}

By formula \eqref{css} and the definition of the function $\phi_\la$, the
function $W_\gamma \cdot \phi_\la$ equals $q^{-2(\la,\rho)} \on{Tr}(\gamma,
V(\la))$ times the characteristic function of the double coset $N \la(\pi)
K$. Hence $s_\gamma(\phi_\la)=\int_{N \backslash G} W_\gamma(g)\phi_\la(g)
dg$ equals $q^{-(\la,\rho)} \on{Tr}(\gamma, V(\la))$ times the measure of
the right $K$--orbit $\la(\pi) \cdot K$ in $N \backslash G$. This measure
equals $\mu(K/Ad_{\la(\pi)}(N(\OO))) =
\mu(N(\OO))/\mu(Ad_{\la(\pi)}(N(\OO)))$ due to our normalization. The
latter equals $q^{2(\la,\rho)}$. Therefore $s_{\gamma}(\phi_\la) =
\on{Tr}(\gamma, V(\la))$ for each $\la\in \PL^+$.

Any $\phi \in \Ql(G/K)^N_\Psi$ can be written as a finite sum $\sum_{\la
\in \PL^+} a_\la \phi_\la$, where $a_\la \in \Ql$. We can identify the
vector space $\Ql(G/K)^N_\Psi$ with $\on{Rep} \GL(\Ql)$, by mapping
$\phi_\la$ to $[V(\la)]$. Let $\ol{\phi}$ be the image of $\phi$ in
$\on{Rep} \GL(\Ql)$ under this identification. Then $s_\gamma(\phi) =
\sum_{\la \in \PL^+} a_\la \on{Tr}(\gamma,V(\la))$ is simply the value of
$\ol{\phi}$ at $\gamma \in \on{Spec} \on{Rep} \GL(\Ql)$. Since the algebra
$\on{Rep} \GL(\Ql)$ has no nilpotents, $\phi=\phi'$, if and only if
$s_\gamma(\phi) = s_\gamma(\phi')$ for all semi-simple conjugacy classes
$\ga$ in $\GL(\Ql)$,

Therefore to prove \thmref{local}, it is sufficient to check that for each
semi-simple conjugacy class $\ga$ in $\GL(\Ql)$ and $\la\in \PL^+$,
$$s_{\gamma}\circ \Phi(H_\la)=\on{Tr}(\gamma,V(\la)).$$ But the
composition $s_{\gamma}\circ \Phi: \Ql(G/K)^K \to \ol{\mathbb
Q}_{\ell,\ga}$ is a homomorphism of ${\mathcal H}$-modules, by
\lemref{F} and \lemref{mapp}. It is easy to check directly that the
value of this homomorphism on the element $H_0 = \on{ch}_K \in
\Ql(G/K)^K$ equals $1$. Therefore
$$s_{\gamma}\circ \Phi(H_\la) = s_\gamma \circ \Phi(H_\la \cdot H_0) =
\chi_{\gamma}(H_\la) \cdot s_{\gamma} \circ \Phi(H_0) =
\chi_{\gamma}(H_\la) = \on{Tr}(\gamma,V(\la))$$ (see \remref{character}),
and \thmref{local} follows.

\begin{rem}    \label{equi}
Our proof shows that \thmref{local} is equivalent to \thmref{prime}.\qed
\end{rem}

\subsection{Proof of \propref{Fplus}}

Note that $(-1)^{2(\la,\rho)} = (-1)^{|\la|(n-1)}$. Hence we obtain from
\propref{hla} and formula \eqref{le1}:
\begin{equation}    \label{lenew} L_{E,x} = \sum_{\la \in P^{++}_n}
q_x^{(\la,\rho)+n(\la)}
\on{Tr}(\sigma(\on{Fr}_x),E_x(\la)) \cdot H_\la.
\end{equation} Therefore we find: $$\int_{N^J(\K_x)} L_{E,x}(u_x \cdot
\on{diag}(\pi_x^{\nu_1},\ldots,\pi_x^{\nu_n})) \Psi^{-1}(u_x) \; du_x =$$
$$\sum_{\la \in P^{++}_n} q_x^{(\la,\rho)+n(\la)}
\on{Tr}(\sigma(\on{Fr}_x),E_x(\la)) \cdot \int_{N^J(\K_x)} H_\la(u_x \cdot
\on{diag}(\pi_x^{\nu_1},\ldots,\pi_x^{\nu_n})) \Psi^{-1}(u_x) \; dx.$$
According to \eqref{formula}, the latter sum equals $q_x^{n(\nu)}
\on{Tr}(\sigma(\on{Fr}_x),E_x(\nu))$, which is the right hand side of
formula \eqref{oh}. Now \propref{Fplus} is proved, and this finishes the
proof of \propref{r} and \thmref{prin}.

\section{Whittaker functions and spherical functions}

In this section we give an interpretation of \thmref{local} from the
point of view of the theory of spherical functions. Throughout this
section we will work over the field of complex numbers instead of
$\Ql$. In particular, all functions will be $\C$--valued, and
${\mathcal H}$ will be a $\C$--algebra.

\subsection{The map $\Th$} Denote by $C^\infty(G/K)^K$ (resp.,
$C^\infty(G/K)^N_\Psi$) the space of smooth left $K$--invariant
(resp., $(N,\Psi)$--equivariant) and right $K$--invariant functions on
$G$. We also denote by $\C(G/K)^K$ (resp., $\C(G/K)^N_\Psi$) the
subspace of compactly supported (resp., compactly supported modulo
$N$) functions.

Each element of $C^\infty(G/K)^K$ can be written as an infinite sum
$\sum_{\la \in \PL^+} a_\la c_\la$, where $c_\la$ is the characteristic
function of the $G(\OO)$--orbit $Gr(G)^\la$.

\begin{lem}    \label{fini}
For each $g \in G$, $(\Phi(c_\la))(g) = 0$ for all but finitely many
$\la \in \PL^+$.
\end{lem}

\begin{proof} It suffices to prove the statement for $g=\mu(\pi)$. In
this case, it is easy to see that for all but finitely many $\la$,
there exists an element $v \in N$ (depending on $\la$) with $\Psi(v)
\neq 1$, such that $\forall u \in N$, $u \cdot \mu(\pi) \in Gr(G)^\la$
if and only if $(vu) \cdot \mu(\pi) \in Gr(G)^\la$. But then
$(\Phi(c_\la))(\mu(\pi)) = \Psi(v) (\Phi(c_\la))(\mu(\pi))$, and
hence $(\Phi(c_\la))(g) = 0$.
\end{proof}

Therefore $\Phi$ defines a map $C^\infty(G/K)^K \arr
C^\infty(G/K)^N_\Psi$, $f \arr \Phi(f)$, which is equivariant with
respect to the action of Hecke operators.

Now we define a map $\Th: C^\infty(G/K)^N_\Psi \arr C^\infty(G/K)^K$ by
the formula
\begin{equation}    \label{inverse}
(\Th(f))(g) = \int_K f(kg) dk,
\end{equation}
where $dk$ stands for the Haar measure on $K$ of volume $1$. This map
is also equivariant with respect to the action of Hecke operators.

We define $a$ as the element of $\C(G/K)^K$ equal to $(\Th \circ
\Phi)(\on{ch}_K) = \Th(\phi_0)$. The same argument as in the proof of
\lemref{fini} shows that the map $\Th$ sends functions from
$\C(G/K)^N_\Psi$ to $\C(G/K)^K$. Hence $a \in \C(G/K)^K = {\mathcal
H}$.

Introduce the notation $$(a*f)(g) = \int_G a(x) f(gx) \; dx.$$
Then we obtain:
\begin{equation}    \label{star}
(\Th \circ \Phi)(f) = a*f, \quad \quad \forall f \in C^\infty(G/K)^K.
\end{equation}

In the next section we will use the element $a$ to clarify the
connection between Whittaker functions and spherical functions.

\subsection{Connection between $a$ and the Plancherel measure}
Let $\gamma$ be a semi-simple conjugacy class in the group
$\GL(\C)$. Recall \cite{Sa,Mac} that the spherical function $S_\gamma$
is the unique $K$--bi-invariant function on $G$, such that

\begin{itemize}
\item $f*S_\ga = \chi_\gamma(f) S_\ga, \forall f \in {\mathcal H}$,
where $\chi_\ga: {\mathcal H} \arr \C$ is the character corresponding
to $\ga$ defined in \remref{character};
\item $S_\ga(1) = 1$.
\end{itemize}

These properties imply that
\begin{equation}    \label{vazh}
\int_G f(x) S_\ga(x) \; dx = \chi_\gamma(f).
\end{equation}

Now let $W_\ga$ be the Whittaker function on $G$ as defined in 
Sect.~5.4 but with the character $\Psi^{-1}$ of $N$ replaced with
the character $\Psi$. It is straightforward to check that the function
$\Phi(S_\gamma)$ satisfies all the properties of the function $W_\ga$
from \secref{proof}, except for the normalization condition
$W_\ga(1)=1$. By \thmref{prime}, $\Phi(S_\gamma)$ is proportional to
$W_\gamma$.

\begin{lem}    \label{whereitgoes}
$$\Phi(S_\ga) = \chi_\ga(a) W_\ga.$$
\end{lem}

\begin{proof}
Introduce $a(\ga)$ by the formula $\Phi(S_\ga) = a(\ga) W_\ga$. Since
$\Th(W_\ga) = S_\ga$ by definition, we obtain, using formula
\eqref{star} and the properties of $S_\ga$: $a(\ga) S_\ga = (\Th \circ
\Phi)(S_\ga) = a*S_\ga = \chi_\ga(a) S_\ga$.
\end{proof}

According to \cite{Mac}, (1.5.1), there exists a unimodular measure
$d\mu(\gamma)$ (Plancherel measure) on the maximal compact subtorus $\TL^u$
of $\TL$, which satisfies
\begin{equation}    \label{plancherel}
\int_G f_1(g) \overline{f_2(g)} dg = \int_{\TL^u} \chi_\ga(f_1)
\overline{\chi_\ga(f_2)} \; d\mu(\gamma),
\end{equation} for all $f_1, f_2 \in {\mathcal H}$.

Setting $f_2=\on{ch}_K$, we obtain:
\begin{equation}    \label{raz}
f(1) = \int_{\TL^u} \chi_\ga(f) \; d\mu(\gamma), \quad \quad \forall f
\in {\mathcal H}.
\end{equation}

By \thmref{local}, $\Phi(H_\la) = \phi_\la$. But it is clear that
$(\Th(\phi_\la))(1) = \delta_{\la,0}$. Therefore, using \eqref{star},
we see that $(a*H_\la)(1)=\delta_{\la,0}$. Substituting this into
formula \eqref{raz} and using the formula $\chi_\ga(H_\la) =
\on{Tr}(\ga,V(\la))$, we obtain:
\begin{equation}    \label{del}
\int_{\TL^u} \on{Tr}(\gamma,V(\la)) a(\gamma) \; d\mu(\gamma) =
\delta_{\la,0}.
\end{equation}

There exists a unique measure $d\wt{\mu}(\gamma)$ on $\TL^u$ (induced
by the Haar measure on $\GL^u$), such that
\begin{equation}    \label{orth}
\int_{\TL^u} \on{Tr}(\gamma,V(\nu)) \ol{\on{Tr}(\gamma,V(\la))} \;
d\wt{\mu}(\gamma) = \delta_{\la,\nu}.
\end{equation} Formula \eqref{del} then implies

\begin{prop}    \label{aga}
$$a(\ga) = \frac{d\wt{\mu}(\gamma)}{d\mu(\gamma)}.$$
\end{prop}

Now \lemref{whereitgoes} gives us:
\begin{equation}    \label{proportion}
\Phi(S_\gamma) = \frac{d\wt{\mu}(\gamma)}{d\mu(\gamma)} W_\gamma.
\end{equation}

\subsection{Another proof of \thmref{local}}
In this subsection, which is independent from the previous one, we use
spherical functions to give another proof of \thmref{local}.

Substituting $f_1=\on{ch}_{KgK}$ into formula \eqref{plancherel} and
using formula \eqref{vazh}, we obtain that for any $f \in {\mathcal
H}$,
\begin{equation}    \label{fourier}
\ol{f(g)} = \int_{\TL^u} S_{\gamma}(g) \ol{\chi_\gamma(f)} \; d\mu(\gamma).
\end{equation}

Since $\chi_\ga(H_\la) = \on{Tr}(\gamma,V(\la))$, we have:
\begin{equation}    \label{alambda} H_\la(g) = \int_{\TL^u} S_\ga(g)
\ol{\on{Tr}(\gamma,V(\la))} \; d\mu(\gamma).
\end{equation}

According to formula \eqref{alambda},
\begin{equation}    \label{a0} H_0(g) = \int_{\TL^u} S_\ga(g) \;
d\mu(\gamma).
\end{equation} Hence $$\Phi(H_0)(g) = \int_{\TL^u} \Phi(S_\ga)(g) \;
d\mu(\gamma) = \int_{\TL^u} W_\ga(g) a(\gamma) \; d\mu(\gamma).$$ On
the other hand, it is clear from definition that $\Phi(H_0) =
\phi_0$. Therefore, substituting $g=\la(\pi)$ and using formula
\eqref{css}, we obtain formula \eqref{del}. Repeating the argument
with the Haar measure given above, we obtain \eqref{proportion}.

Now formulas \eqref{alambda}, \eqref{proportion} and \eqref{css} give:
$$(\Phi(H_\la))(\nu(\pi)) = \int_{\TL^u} W_\gamma(\nu(\pi))
\ol{\on{Tr}(\gamma,V(\la))} \; d\wt{\mu}(\gamma) =$$ $$=
q^{-(\nu,\rho)} \int_{\TL^u} \on{Tr}(\gamma,V(\nu))
\ol{\on{Tr}(\gamma,V(\la))} \; d\wt{\mu}(\gamma) = q^{-(\la,\rho)}
\delta_{\la,\nu}.$$ This proves formula \eqref{formula} and
\thmref{local} over the field of complex numbers. Since $H_\la$ takes
values in rational numbers and $\Psi$ takes values in the roots of
unity, the validity of \eqref{formula} over $\C$ actually implies its
validity over $\Ql$.

\subsection{The function $L_\ga$}
The Whittaker function can be written as a series
$$W_\ga = \sum_{\la \in \PL^+} \on{Tr}(\gamma,V(\la)) \cdot
\phi_\la.$$ This series obviously makes sense, since the supports of
the functions $\phi_\la$ do not intersect. In view of \thmref{local},
it is natural to consider the series
\begin{equation}    \label{Fgamma}
L_\ga = \sum_{\la \in \PL^+} \on{Tr}(\gamma,V(\la)) \cdot H_\la.
\end{equation}
However, the convergence of this series is not at all automatic,
because the supports of functions $H_\la$ do intersect; for instance,
each $H_\la$ has a non-zero value at $1$. In this section we study the
question of convergence of $L_\ga$.

Let us write: $$H_\la = q^{-(\la,\rho)} \sum_{\mu \leq \la}
P_{\mu\la}(q^{-1}) \cdot c_\mu,$$ where $q^{-(\la,\rho)} P_{\mu\la}$ is a
polynomial in $q^{-1}$ (recall that $c_\la$ is the characteristic function
of the $G(\OO)$--orbit $\G(G)^\la$). Formula \eqref{alambda} can be
rewritten as follows:
\begin{equation}    \label{al1}
H_\la(g) = \int_{\TL^u} S_\ga(g) a(\ga)^{-1}
\ol{\on{Tr}(\gamma,V(\la))} \; d\wt{\mu}(\gamma).
\end{equation}
Using the defining properties of the spherical function $S_\ga$, we
can write it as a series
\begin{equation}    \label{sga}
S_\ga = \sum_{\mu \in \PL^+} s_\ga^\mu(q^{-1}) \cdot c_\mu.
\end{equation}
where $s_\ga^\mu(q^{-1})$ is a rational function in $q^{-1}$ of the
form $Q(q^{-1}) \wt{s}_\ga^\mu(q^{-1})$. Here $$Q(q^{-1}) = \prod_{i=1}^l
\frac{1-q^{-m_i-1}}{1-q^{-1}}$$ ($l$ is the rank of $G$, $m_i$'s are
the exponents of $G$; note that $Q(q) = \# G/B(\Fq)$), and
$\wt{s}_\ga^\mu$ is a polynomial in $q^{\pm 1}$. Its coefficients are
finite integral linear combinations of characters of irreducible
representations of $\GL$ (for an explicit formula, see \cite{Mac}). It
follows from formula \eqref{whereitgoes} that $a(\ga)$ has the same
structure as a function of $q^{-1}$. Hence both $s_\ga^\mu(q^{-1})$
and $a(\ga)^{-1}$ can be viewed as formal Laurent power series in
$q^{-1}$ and formula \eqref{al1} can be viewed as an identity on such
power series.

We have: $$s_\ga^\mu = \sum_{m > -M} \on{Tr}(\ga,R^\mu_m) q^{-m},$$
and $$a(\ga)^{-1} = \sum_{m > -M'} \on{Tr}(\ga,U_m) q^{-m},$$ where
$R^\mu_m$ and $U_m$ are finite linear combinations of irreducible
representations of $\GL$ (the summation is actually only over $m \in
\Z_+$). Then formula \eqref{al1} can be written as follows:
\begin{equation}    \label{f1}
q^{-(\la,\rho)} P_{\mu\la}(q^{-1}) = \sum_{N \in \Z} q^{-N} \int_{\TL^u}
\on{Tr}(\ga,\oplus_{m \in \Z} R^\mu_m \otimes U_{N-m})
\ol{\on{Tr}(\gamma,V(\la))} \; d\wt{\mu}(\gamma).
\end{equation}
By formula \eqref{orth}, the $q^{-N}$ coefficient of
$P_{\mu\la}(q^{-1})$ equals the multiplicity of $V(\la)$ in $\oplus_{m
\in \Z} R^\mu_m \otimes U_{N-m}$. But by construction the latter is
a finite linear combination of irreducible representations of
$\GL$. Hence we obtain the following

\begin{lem}
For each $\mu \in \PL^+$ and $N \in \Z_+$ there are only finitely many $\la
\in \PL^+$, such that $q^{-(\la,\rho)} P_{\mu\la}$ has a non-zero
coefficient in front of $q^{-N}$.
\end{lem}

\begin{rem} This can also be seen from the explicit formula for
$P_{\mu\la}$ obtained in \cite{Lu2,Ka}.\qed
\end{rem}

Therefore for each $g \in G$, $L_\ga(g)$ given by formula
\eqref{Fgamma} makes sense as a formal power series in $q^{-1}$, each
coefficient being a finite linear combination of
characters. Furthermore, we see, by reversing the argument above that
as formal power series in $q^{-1}$,
\begin{equation}    \label{propor}
L_\ga(g) = a(\ga)^{-1} S_\gamma(g), \quad \quad \forall g \in G.
\end{equation}

In order to estimate  the convergence of this series, we have to compute
$a(\ga)$ explicitly. According to \cite{Mac}, (5.1.2),
$$d\mu(\gamma) = \frac{Q(q^{-1})}{|W|} \frac{\prod_{\al \in \De}
(1-\al(\gamma))}{\prod_{\al \in \De} (1-q^{-1} \al(\gamma))} d\gamma,$$
where $|W|$ is the number of elements in the Weyl group, and $\De$
is the set of roots of $G$. The notation $d\ga$ means the Haar measure on
$\TL^u$, which gives it volume $1$. On the other hand,
$$d\wt{\mu}(\gamma) = \frac{1}{|W|} \prod_{\al \in \De} (1-\al(\gamma))
d\gamma.$$ Now \propref{aga} gives:
\begin{equation}    \label{propsg}
a(\ga) = \frac{\prod_{\al \in \De} (1-q^{-1} \al(\gamma))}{Q(q^{-1})}.
\end{equation}

Let $^L\!{\mathfrak g}$ be the adjoint representation of $\GL$, and
$\Lambda^i(^L\!{\mathfrak g})$ be its $i$th exterior power. Formula
\eqref{propsg} means that $a \in {\mathcal H}$ is the image under the
Satake isomorphism $S$ of the following element in $\on{Rep} \GL$:
$$\prod_{i=1}^l (1-q^{-m_i-1})^{-1} \sum_{i=0}^{\dim ^L\!{\mathfrak
g}} \left[ \Lambda^i(^L\!{\mathfrak g}) \right] (-1)^i q^{-i}.$$

Formula \eqref{propor} implies the following result.

\begin{prop}    \label{converge}
If the conjugacy class $\ga$ satisfies: $q^{-1} < |\al(\ga)| < q,
\forall \al \in \De_+$, then $L_\ga(g)$ converges absolutely to
$$\frac{Q(q^{-1})}{\prod_{\al \in \De} (1-q^{-1} \al(\gamma))}
S_\ga(g)$$ for all $g \in G$.
\end{prop}

Note that when $G=GL_n$, formula \eqref{Fgamma} looks similar to
formula \eqref{lenew} for the function $L_{E,x}$. Besides powers of
$q_x$, the difference is that in  \eqref{lenew} the summation is
restricted to the subset $P^{++}_n$ of the set $P^+_n = \PL^+$ of all
dominant weights of $GL_n$. However, $L_{E,x}$ is not the restriction
of \eqref{Fgamma} to the union of strata $Gr(GL_n)^\la$ with $\la \in
P^{++}_n$, because the functions $H_\la$ with $\la \in P^+_n -
P^{++}_n$ do not vanish on those strata. While $L_\gamma$ given by
\eqref{Fgamma} is manifestly an eigenfunction of the Hecke operators,
$L_{E,x}$ is not\footnote{it is actually an eigenfunction of some
other operators, similar to the Hecke operators, which were defined by
Laumon \cite{La1}}.

For general $G$ there is no analogue of the subset $\PL^{++} \subset
\PL^+$, and so the function $L_\gamma$ seems to be the closest
analogue of $L_{E,x}$ in the general setting. According to
\thmref{local} and formula \eqref{css}, $\Phi(L_\gamma)$ equals the
Whittaker function $W_\ga$.

\begin{rem}    \label{lfunct}
Let $\ga$ be a semi-simple conjugacy class of $\GL(\Ql)$ and $r:
\GL(\Ql) \arr \on{Aut} V$ be a finite-dimensional representation of
$\GL(\Ql)$. Recall that the local $L$--function associated to the pair
$(\ga,V)$ is defined by the formula
\begin{equation}    \label{lfun} L(\ga,V;s) = \on{det} \left( 1 - r(\ga)
q^{-s} \right)^{-1}.
\end{equation} In particular, if $V = ^L\!{\mathfrak g}$ is the adjoint
representation, then
$$L(\ga,^L\!{\mathfrak g};s) = (1-q^{-s})^{-l} \prod_{\al \in \De} (1 -
\al(\ga) q^{-s})^{-1}.$$ Hence $a(\ga)$ can be written as
$$a(\ga) = L(\ga,^L\!{\mathfrak g};1) \prod_{i=1}^l (1-q^{-m_i-1}).$$
Thus, we obtain:
\begin{equation}    \label{vani}
\Phi(S_\ga) = L(\ga,^L\!{\mathfrak g};1)^{-1} \prod_{i=1}^l
(1-q^{-m_i-1})^{-1} \cdot W_\ga.
\end{equation}

Using arguments similar to those of Casselman and Shalika \cite{CS}, one
can show that the irreducible unramified representation corresponding to
$\ga$ has a Whittaker model if and only if $\Phi(S_\ga) \neq 0$. Formula
\eqref{vani} means that $\Phi(S_\ga) \neq 0$ if and only if
$L(\ga,^L\!{\mathfrak g};s)$ is regular at $s=1$. We conclude that the
irreducible unramified representation of $G$ with the Langlands parameter
$\ga$ has a Whittaker model if and only if $L(\ga,^L\!{\mathfrak g};s)$ is
regular at $s=1$ (in that case the Whittaker model is actually
unique). This agrees with a special case of a conjecture of Gross and
Prasad \cite{GP} (Conjecture 2.6).\qed
\end{rem}

\subsection{Identities on $P_{\mu\la}$} According to formula
\eqref{propor}, for each $\mu \in \PL^+$ we have the following
equality of power series in $q^{-1}$:
\begin{equation}    \label{id1}
\sum_{\la: \la \geq \mu} \on{Tr}(\gamma,V(\la)) q^{-(\la,\rho)}
P_{\mu\la}(q^{-1}) = \prod_{i=1}^l \frac{1-q^{-m_i-1}}{1-q^{-1}} \prod_{\al
\in \De} (1-q^{-1} \al(\gamma))^{-1} S_\ga(\mu(\pi)).
\end{equation}

Recall that the coefficients of the polynomial $P_{\mu\la}$ (which can
be interpreted as a Kazhdan-Lusztig polynomial for the affine Weyl
group \cite{Lu2,Ka}) are given by dimensions of stalk cohomologies of
the perverse sheaf ${\mathcal A}_\la$:
\begin{equation}
P_{\mu\la} = q^{(\la,\rho)} \sum_{i \in \Z} \dim \on{H}^i({\mathcal
A}_\la)|_{\mu(\pi)} \; q^{i/2}.
\end{equation}
Thus, formula \eqref{id1} is an identity which connects these
dimensions with the values of the spherical functions. The latter are
known explicitly; they can be expressed via the Hall-Littlewood
polynomials \cite{Mac}.

For example, let us apply formula \eqref{id1} when $G=SL_2$ and
$\mu=0$. In this case, $\gamma \in \C^\times$, and the set $\PL^+$ of
dominant weights of the dual group $\GL=PGL_2$ can be identified with
the set of non-negative even integers. To weight $2m$ corresponds the
$2m+1$--dimensional representation $V(2m)$ of $PGL_2$, and
$\on{Tr}(\gamma,V(2m)) = \sum_{i=-m}^m \gamma^i$. Formula \eqref{id1}
then gives:
$$\sum_{m \in \Z_+} \left( \sum_{j=-m}^m \gamma^j \right) \cdot
\sum_{i \in \Z} \dim \on{H}^i({\mathcal A}_{2m})|_1 \; q^{i/2} =
\frac{1+q^{-1}}{(1-q^{-1} \gamma)(1-q^{-1} \gamma^{-1})}.$$ This is
easy to see directly, because ${\mathcal A}_{2m}$ is known to be the
constant sheaf on the closure of the stratum $\G(SL_2)^{2m}$ placed in
degree $-2m$.

For general $G$, formula \eqref{id1} with $\mu=0$ can be interpreted as
follows. Let $R(\g)$ be the graded ring of polynomials on $\g$, $J(\g)$ be
its subgring of $\GL$--invariants, and $H(\g)$ be the subspace of
$\GL$--harmonic polynomials on $\g$.

For a graded space $V$, we denote by $V_j$ its $j$th homogeneous
component. If each $V_j$ is a representation of $\GL$, we denote by
$\on{Ch}(\ga,V)$ the graded character of $V$:
$$\on{Ch}(\ga,V) = \sum_{j=0}^\infty \on{Tr}(\ga,V_j) q^{-j}.$$
Note that the graded character of $R(\g)$ equals
\begin{equation}    \label{tozh1}
\on{Ch}(\ga,R(\g)) = \prod_{i=1}^l (1-q^{-1})^{-l} \prod_{\al \in \De}
(1-q^{-1} \al(\gamma))^{-1}
\end{equation}
(compare with \remref{lfunct}), while
$$\on{Ch}(\ga,J(\g)) = \prod_{i=1}^l (1-q^{-m_i-1})^{-1}.$$

Now let $H(\g)$ be the (graded) space of $\GL$--harmonic polynomials on
$\g$. According to Theorem 0.2 of B.~Kostant \cite{Ko}, $R(\g) = J(\g)
\otimes H(\g)$. Hence
$$\on{Ch}(\ga,H(\g)) = \prod_{i=1}^l \frac{1-q^{-m_i-1}}{1-q^{-1}}
\prod_{\al \in \De} (1-q^{-1} \al(\gamma))^{-1} = a(\ga)^{-1}.$$ Thus, we
obtain another interpretation of $a(\ga)^{-1}$: it is equal to the graded
character of the space of $\GL$--harmonic polynomials. Note that  it also
equals the graded character of the space of regular functions on the
nilpotent cone ${\mathcal N}$ in $\g$.

Formula \eqref{id1} together with this interpretation give us the following
result.

\begin{prop}
\begin{equation}    \label{new}
P_{0,\la}(q^{-1}) = q^{(\la,\rho)} \sum_{j=0}^\infty q^{-j}
\on{mult}(V(\la),H(\g)_j),
\end{equation}
where $\on{mult}(V(\la),H(\g)_j)$ is the multiplicity of $V(\la)$ in
 $H(\g)_j$.
\end{prop}

A complete description of these multiplicities has been given by Kostant in
\cite{Ko}. In fact, applying Theorem 0.11 of \cite{Ko} to formula
\eqref{new}, we obtain:
\begin{equation}    \label{tozh2}
P_{0,\la} = q^{(\la,\rho)} \sum_{i=1}^{l_\la} q^{-m_i(\la)},
\end{equation}
where $m_i(\la)$ are the generalized exponents associated to the
representation $V(\la)$, defined in \cite{Ko}. In the special case when
$V(\la)$ is the adjoint representation $^L\!{\mathfrak g}$, these are just
the exponents of $\GL$, and formula \eqref{tozh2} specializes to Lusztig's
formula (see \cite{Lu2}, p.~226)
$$P_{0,\la_{\on{adj}}} = \sum_{i=1}^l q^{m_i-1}.$$

Formula \eqref{tozh2} is not new: R.~Brylinski \cite{Br} observed that
it immediately follows if one compares the Lusztig-Kato formula
\cite{Lu2,Ka} for $P_{0,\la}$ and the Hesselink-Peterson formula
\cite{He} for the right hand side of \eqref{tozh2}. Note that in
contrast to her argument, our proof of formula \eqref{tozh2} does
not use  the Lusztig-Kato formula.

Using \eqref{f1} it is easy to derive a formula analogous to \eqref{new}
for $P_{\mu\la}$ in terms of Hall-Littlewood polynomials.

\section{Geometric analogue of \thmref{local} and some open problems}

\subsection{}
In this subsection we will formulate a geometric analogue of
\thmref{local}. Recall the Grassmannian $\G(G)$ of section ~5.2. This
is a strict ind--scheme over $\Fq$, i.e., an inductive system of
$\Fq$--schemes $\G(G)_k, k\geq 0$, where all maps $i_{k,m}: \G(G)_k
\rightarrow \G(G)_m, k<m$, are closed embeddings. For more details,
see, e.g., \cite{LS,MV}. By a $\Ql$--sheaf on $\G(G)$ we will
understand a system of $\Ql$--sheaves $\cf_k$ on each $\G(G)_k$ and a
compatible system of isomorphisms $\cf_k \simeq i_{k,m}^* \cf_m$ for
 all $k<m$.

There exists an ind--group scheme ${\mathcal N}(\K)$, whose set of
$\Fq$--points is $N(\K)$. This ind-group scheme acts on the Grassmannian
$\G(G)$, and its orbits stratify $\G(G)$ by ind-schemes $S^\nu, \nu \in
\PL$. The stratum $S^\nu$ is the $\NK$--orbit of $\nu(\pi) \in
Gr(G)$. Denote by $j^\nu$ the embedding $S^\nu \hookrightarrow \G(G)$.

We choose a generic additive character
$\wt\Psi: \NK \to {\mathbb G}_{a,\mathcal K}$ and define a
homomorphism $\Psi: \NK \to \Gaf$ by the formula
$$
\Psi \ = \ \operatorname{Res} \circ \,\wt\Psi: \NK \to \Gaf\,,
$$
where $\operatorname{Res}$ is the geometric analogue of the residue
map of Sect.~5.3. As before, let $\psi: \Fq \to \Ql^\times$ denote a
non-trivial character and let $\mathcal I_\psi$ denote the
corresponding Artin-Schreier sheaf on $\Gaf$.

Consider the category $\on{P}_{\GG(\OO)}(\G(G))$ of
${\GG(\OO)}$--equivariant perverse sheaves on $\G(G)$ with
finite-dimensional support. This category is a geometric analogue of
the Hecke algebra ${\mathcal H}=\Ql(G/K)^K$ (see \remref{cate}
below). We define an abelian category $\on{Sh}^\Psi_{\NK}(\G(G))$ (a
geometric analogue of $\Ql(G/K)^N_\Psi$) and a collection of
cohomology functors $\W^i: \on{P}_{\GG(\OO)}(\G(G)) \arr
\on{Sh}^\Psi_{\NK}(\G(G))$, which are a geometric analogue of the map
$\W$ of Sect.~5.3.

The objects of the category $\on{Sh}^\Psi_{\NK}(\G(G))$ are
$\Ql$--sheaves ${\mathcal E}$ on $\G(G)$ which satisfy the following
conditions:

(1) ${j^\nu}^* {\mathcal E} = 0$ except for finitely many $\nu\in \PL$;

(2) $t_\nu^* {j^\nu}^* {\mathcal E} \otimes\varphi^*{\mathcal
I}_{\psi^{-1}}$ are trivial local systems of finite rank for all $\nu$,
where $t_\nu: \NK \to S^\nu$ is the map $u \mapsto u \cdot\nu(t)$.

Morphisms in this category are defined in an obvious way.

\begin{lem} \label{first}
If $\nu\in\PL$ is not dominant, then for every ${\mathcal E} \in\on{Ob}
(\on{Sh}^\Psi_{\NK}\G(G))$, ${j^\nu}^*({\mathcal E})=0$.
\end{lem}

The proof is analogous to the proof of the corresponding statement for
functions. Thus, ${\mathcal E}$ not only satisfies property (1) above, but
also satisfies the stronger property

(1') ${j^\nu}^* {\mathcal E} = 0$ except for finitely many $\nu\in \PL^+$.

Now we construct the functors  $\W^i:\on{P}_{\GG(\OO)}(\G(G))\to
\on{Sh}^\Psi_{\NK}(\G(G))$. Consider the sequence of maps:
$$
\begin{CD}
\G @<a<< \NK\times\G @>q>> \NK/\NO\times\G @>p>> \G\,,
\end{CD}
$$ where $a$ is given by acting with $\NK$ on $\G$ and $p,q$ are
projections. For $\cf\in \on{P}_{\GG(\OO)}(\G(G))$ we then set:
$${\W}^i(\cf) = R^i p_!(\wt \cf \otimes \varphi^*{\mathcal I}_\psi), \qquad
\text{with} \ \ q^*\wt\cf = a^*\cf\,.
$$

To formulate the geometric analogue of \thmref{local}, recall that for
$\la\in {^LP^+}$ we denote by ${\mathcal A}_\lambda$ the Goresky-MacPherson
extension of the sheaf $\ol{\mathbb Q}_{\ell,\la}[2(\la,\rho)]((\la,\rho))$
associated to the $\GG(\OO)$--orbit $\G(G)^\la$. For each $\nu \in \PL^+$,
denote by $\NK_\nu$ the isotropy group of $\nu(\pi)$. Since $\nu$ is
dominant, the restriction of $\Psi$ to $\NK_\nu$ equals $0$. Therefore the
map $\Psi$ restricts to a map on $S^\nu$, which we continue to denote by
the same letter $\Psi: S^\nu\to\Gaf$. With this notation, the sheaf
$\Psi^*{\mathcal I}_\psi$ is a sheaf on $S^\nu$. Recall that $j^\nu$
denotes the embedding $S^\nu \hookrightarrow \G(G)$. Now we are ready to
state the geometric analogue of \thmref{local}.

\begin{conj} \label{general2}
$${\W}^i({\mathcal A}_\lambda) = \left\{\aligned &j^\la_!  \Psi^*{\mathcal
I}_\psi((\la,\rho)) \qquad\text{if} \ \ i=2(\lambda,\rho) \\ &\ 0
\qquad\text{if} \ \ i\neq 2(\lambda,\rho) \,.\endaligned\right.
$$
\end{conj}

We will now formulate \conjref{second} describing the stalk
cohomologies of the sheaves ${\W}^i({\mathcal A}_\lambda)$. The
statement of \conjref{second} does not explicitly involve the category
$\on{Sh}^\Psi_{\NK}(\G(G))$ and the functors $\Phi^i$. However,
\conjref{general2} can be derived from \conjref{second}.

Note that the support of the restriction of ${\mathcal A}_\la$ to $S^\nu$,
i.e., $\ol{\G(G)^\la} \cap S^\nu$, is finite-dimensional. To work in a
geometric setting, let us extend the base field from $\Fq$ to $\ol{\mathbb
F}_q$ and use Weil sheaves. Denote by $\Psi_{\la,\nu}$ the restriction of
$\Psi$ to $\ol{\G(G)^\la} \cap S^\nu$.

\begin{conj}    \label{second} For $\lambda$ dominant
$$\operatorname{H}_c^k(\ol{\G(G)^\la} \cap S^\nu, {\mathcal
A}_\lambda\otimes \Psi_{\la,\nu}^*{\mathcal I}_\psi) =
\left\{\aligned &\Ql(-(\la,\rho)) \ \ \text{if} \ \nu=\lambda \
\text{and} \ \ k=2(\la,\rho)
\\ &\ 0\qquad\qquad\ \ \text{otherwise}\,.\endaligned \right.
$$
\end{conj}

Proving \conjref{second} would yield an alternative proof of
\thmref{local}, and hence of the Casselman-Shalika formula \eqref{css}
(see \remref{equi}).

One sees readily that \conjref{second} holds when $\la \leq \nu$. Let us,
then, assume that $\nu$ is dominant and $\nu<\lambda$. The projection
formula implies that
$$
\operatorname {H}_c^*(S^\nu, {\mathcal A}_\lambda\otimes \Psi^*{\mathcal
I}_\psi) = \operatorname {H}_c^*(\Gafb, R\Psi_!{\mathcal A}_\lambda\otimes
{\mathcal I}_\psi).
$$
Therefore \conjref{second} follows from the following

\begin{conj}    \label{third} For $\lambda, \nu$ dominant and
$\nu\neq\lambda$ the sheaves $ R^k\Psi_!{\mathcal A}_\lambda$ are
constant on $\Gafb$.
\end{conj}

By theorem 4.3a of \cite{MV} we see that if \conjref{third} holds then
$R^k\Psi_!{\mathcal A}_\lambda = 0$ unless $k = 2((\la,\rho)-1)$. Here we
are using the fact that the results of \cite{MV}, stated there over
$\mathbb C$, extend to our current context.

\subsection{}
In this subsection we speculate about what could be the analogue of
Laumon's sheaf ${\mathcal L}_E$ in the case of an arbitrary reductive
group. Recall from Sect.~3.2 that ${\mathcal L}_E$ is a sheaf on the
stack $\Coh_n$ canonically attached to a rank $n$ local system $E$ on
$X$. This sheaf is used as the starting point of the conjectural
construction of the automorphic sheaf on ${\mathcal M}_n$ associated
to $E$, see \cite{La2} and Sect.~3 above.

First we revisit the case of $GL_n$ and introduce a scheme
$\G^+_{X^{(\infty)}}$ with a smooth morphism $q: \G^+_{X^{(\infty)}}
\arr \Coh_n$, and take the pull-back of ${\mathcal L}_E$ to
$\G^+_{X^{(\infty)}}$. The scheme $\G^+_{X^{(\infty)}}$ classifies
pairs $\{ L,t \}$, where $L$ is a rank $n$ bundle on $X$ and $t:
\OO_X^{\oplus n} \arr L$ is an embedding of $\OO_X$--modules. The
scheme $\G^+_{X^{(\infty)}}$ is a disjoint union of the smooth schemes
$\G^{+,m}_{X^{(\infty)}}, m\geq 0$, corresponding to bundles of degree
$m$. The scheme $\G^{+,m}_{X^{(\infty)}}$ is isomorphic to the
Grothendieck Quot--scheme $\on{Quot}^m_{\OO^{\oplus
n}_X/X/\Fq}$. Recall \cite{Gr} that $\on{Quot}^m_{\OO^{\oplus
n}_X/X/k}$ classifies the quotients of $\OO^{\oplus n}_X$ that are
torsion sheaves of length $m$; at the level of points, $\{ L,t \}$
corresponds to the quotient of $\OO^{\oplus n}_X$ by the image of
$L^*$ under the transpose homomorphism $t^*: L^* \arr \OO^{\oplus
n}_X$.

The morphism $q: \G^+_{X^{(\infty)}} \arr \Coh_n$ sends $\{ L,t \}$ to
$\OO^{\oplus n}_X/\on{Im} t^*$. In the same way as in Sect.~4.2, one
can show that $q$ is smooth. We denote by the same character
${\mathcal L}_E$ the pull-back of ${\mathcal L}_E$ by $q$. It is the
pair $(\G^+_{X^{(\infty)}},{\mathcal L}_E)$ that we would like to
generalize to other groups.

Let us describe the basic structure of $\G^+_{X^{(\infty)}}$. For each
$k\geq 1$, we introduce the scheme $\G^+_{X^{(k)}}$ over
$X^{(k),rss}$ (see Sect.~3.2), which parametrizes the objects $\{
(x_1,\ldots,x_k),$ $L,t \}$, where $(x_1,\ldots,x_k)$ is a set of $k$
non-ordered distinct points of $X$, $L$ is a rank $n$ bundle over $X$, and
$t$ is its trivialization over $X - \{ x_1,\ldots,x_k \}$, which extends to
an embedding of $\OO_X$--modules $\OO_X^{\oplus n} \arr L$. The fiber of
this scheme over $(x_1,\ldots,x_k) \in X^{(k),rss}$ is the product of  the
$\G^+_{x_i}$. It is easy to describe the pull-back
${{\mathcal L}}_{X^{(k)}}^{E,+}$ of ${\mathcal L}_E$ under the natural
morphism
$\delta_k: \G^+_{X^{(k)}} \arr \G^+_{X^{(\infty)}}$. In particular, the
restriction of ${{\mathcal L}}_{X^{(k)}}^{E,+}$ to the fiber of
$\G^+_{X^{(k)}}$ over $(x_1,\ldots,x_k)$ is $\boxtimes_{i=1}^k {{\mathcal
L}}_{x_i}^{E,+}$, where
\begin{equation}    \label{box1} {{\mathcal L}}_x^{E,+} = \sum_{\la \in
P^{++}_n} {\mathcal A}_{\la,x}[|\la|(n-1)](|\la|(n-1)/2) \otimes E_x(\la).
\end{equation}

The set of $\Fq$--points of $\G^+_{X^{(\infty)}}$ is isomorphic to the
quotient $GL_n(\A)^+/GL_n(\OO)$. For groups other than $GL_n$ we do not
have analogues of the subset $GL_n(\A)^+ \subset GL_n(\A)$, the subset
$P^{++}_n \subset \PL^+$, and the subscheme $\G^+$ of the affine
Grassmannian. For this reason, we can not avoid considering a substantially
larger object in place of $\G^+_{X^{(\infty)}}$.

Thus, for a reductive group $G$, we consider the set
$G(\A)/G(\OO)$. This set, which we denote by $Gr(G)_{X^{(\infty)}}$,
is isomorphic to the set of isomorphism classes of pairs $\{ {\mathcal
P},t \}$, where ${\mathcal P}$ is a principal $G$--bundle over $X$ and
$t$ is an isomorphism between ${\mathcal P}$ and the trivial bundle on
a Zariski open subset of $X$. It is not difficult to define a functor
$\G(G)_{X^{(\infty)}}$ from the category of $\Fq$--schemes to the
category of sets, whose set of $\Fq$--points is
$Gr(G)_{X^{(\infty)}}$. It would be desirable to have a notion of
perverse sheaf on $\G(G)_{X^{(\infty)}}$. A $\GL$--local system $E$ on
$X$ should give rise to a perverse sheaf ${\mathcal L}_E$ on
$\G(G)_{X^{(\infty)}}$, analogous to the sheaf ${\mathcal L}_E$ in the
case of $GL_n$; this sheaf should be irreducible if $E$ is
irreducible. Although we do not know how to define such an object, we
describe below what its pull-backs should be under certain natural
morphisms.

For each $k\geq 1$, following Beilinson and Drinfeld, we introduce the
ind--scheme $\G(G)_{X^{(k)}}$ over $X^{(k),rss}$, which parametrizes
the objects $\{ (x_1,\ldots,x_k),{\mathcal P},t \}$, where
$(x_1,\ldots,x_k)$ is a set of non-ordered distinct points of $X$,
${\mathcal P}$ is a principal $G$--bundle over $X$, and $t$ is its
isomorphism with the trivial bundle over $X - \{ x_1,\ldots,x_k
\}$. The fiber of this scheme over $(x_1,\ldots,x_k) \in X^{(k),rss}$
is the product of  the $\G_{x_i}$ (see \cite{MV}, Sect.~3). We have an
obvious set-theoretic map $\delta_k: Gr(G)_{X^{(k)}} \arr
G(G)_{X^{(\infty)}}$, which can also be defined on the level of
functors: schemes $\arr$ sets. The pull-back of ${\mathcal L}_E$ to
$\G(G)_{X^{(k)}}$ should be the sheaf ${{\mathcal L}}_{X^{(k)}}^E$ on
$\G(G)_{X^{(k)}}$ (inductive limit of perverse sheaves), such that its
restriction to the fiber $\prod_{i=1}^k \G(G)_{x_i}$ over
$(x_1,\ldots,x_k)$ is $\boxtimes_{i=1}^k {{\mathcal L}}_{x_i}^E$,
where
\begin{equation}    \label{rx} {{\mathcal L}}_x^E = \oplus_{\la \in \PL^+}
{\mathcal A}_{\la,x} \otimes E_x(\la)^*
\end{equation} (up to shifts in degree). Here $E_x(\la)$ has the same
meaning as in the case of $GL_n$.

The sheaves ${{\mathcal L}}_{X^{(k)}}^E$ have been previously
considered by Beilinson and Drinfeld in the context of the geometric
Langlands correspondence. Formula \eqref{rx} is analogous to formula
\eqref{box1}. The main difference is that in \eqref{box1} the
summation is restricted to the subset $P^{++}_n$ of the set $P^+_n =
\PL^+$ of all dominant weights of $GL_n$, which does not have an
analogue for general $G$ (compare with Sect.~6.4).

\begin{rem}    \label{cate}
Let $\on{P}_{\GG(\OO_x)}(\G(G)_x)$ be the category of
$\GG(\OO_x)$--equivariant perverse sheaves on $\G(G)_x$ (we consider
objects of $\on{P}_{\GG(\OO_x)}(\G(G)_x)$ as pure of weight $0$). This
category is a tensor category, and as such, it is equivalent to the
tensor category ${{\mathcal R}ep} \GL$ of finite-dimensional
representations of $\GL$. To be precise, this result has been proved
in \cite{Gi,MV} over the ground field $\C$ (in this setting, this
isomorphism was conjectured by V.~Drinfeld; see also \cite{Lu2}). But
the proof outlined in \cite{MV} can be generalized to the $\Fq$--case,
so here we assume the result to be true over $\Fq$ as well.

Note that there is a small error in \cite{MV}. The tensor structure
(or, more precisely, the commutativity constraint), which is given by
the convolution product, should be altered slightly. This alteration
does not affect the structure of $\on{P}_{\GG(\OO_x)}(\G(G)_x)$ as a
monoidal category. We simply replace the perverse sheaf ${\mathcal
A}_{\la,x}$ with $(-1)^{2(\la,\rho)}{\mathcal A}_{\la,x}$, where the
sign $(-1)^{2(\la,\rho)}$ is to be viewed as a formal symbol. The
symbol $(-1)^{2(\la,\rho)}$ has the effect of making the cohomology of
$(-1)^{2(\la,\rho)}{\mathcal A}_{\la,x}$ lie in even degrees
only. Then the equivalence of categories ${{\mathcal R}ep} \GL \arr
\on{P}_{\GG(\OO_x)}(\G(G)_x)$ takes the irreducible representation
$V(\la)$ to the perverse sheaf $(-1)^{2(\la,\rho)}{\mathcal
A}_{\la,x}$. With this adjustment the sign in \propref{hla}
disappears and the equivalence above can be viewed as a categorical
version of the Satake isomorphism $\on{Rep} \GL \arr {\mathcal H}$
(see \thmref{Satake}). Indeed, an equivalence of two categories
induces an isomorphism of their Grothendieck rings. But the
Grothendieck ring of $\on{P}_{\GG(\OO_x)}(\G(G)_x)$ is canonically
isomorphic to the Hecke algebra ${\mathcal H}$ via the
``faisceaux--fonctions'' correspondence.

Now consider the left regular representation of $\GL$, $$\oplus_{\la
\in \PL^+} V(\la) \otimes V(\la)^*$$ as an ind--object of the category
${{\mathcal R}ep} \GL$. The corresponding ind--object of the category
$\on{P}_{\GG(\OO_x)}(\G(G)_x)$ is ${{\mathcal L}}^E_x$ given by
formula \eqref{rx}, adjusted for the formal signs, i.e.,
\begin{equation}    \label{rxx} {{\mathcal L}}_x^E = \oplus_{\la \in \PL^+}
(-1)^{2(\la,\rho)}{\mathcal A}_{\la,x} \otimes E_x(\la)^*.
\end{equation}
\qed
\end{rem}

\begin{rem} Let $L_X^E$ be the function associated to the sheaf ${{\mathcal
L}}_X^E$, and $L_x^E$ be the restriction of $L_X^E$ to $Gr(G)_x \subset
Gr(G)_X$. Using \propref{hla} we obtain:
\begin{equation}    \label{series} L_x^E = \sum_{\la \in \PL^+}
\on{Tr}(\gamma_x,V(\la)^*) \cdot H_{\la,x},
\end{equation} where $\gamma_x = \sigma(\on{Fr}_x)$. Hence the function
$L_x^E$ coincides with the function $L_{\gamma_x^{-1}}$ given by
formula \eqref{Fgamma}. According to \propref{converge}, the series
\eqref{series} converges absolutely, if $q_x^{-1} < |\al(\gamma_x)| <
q_x, \forall \al \in \De_+$, and is proportional to the spherical
function $S_{\gamma_x^{-1}}$.\qed
\end{rem}

\end{document}